%% file: v22.tex
\documentclass[12pt,oneside]{amsart}
\usepackage{mathcomp,longtable, hhline, array}
\usepackage{amsfonts,amssymb,mathrsfs}
\usepackage{url} 
\usepackage{hyperref} 
\usepackage[matrix,arrow,curve]{xy}
\usepackage{upgreek,textcomp}
\usepackage{enumitem}
\setlist[enumerate]{leftmargin=25pt}

\makeatletter 
\@addtoreset{equation}{section} 
\makeatother

\newcommand{\comp}{\mathbin{\scriptstyle{\circ}}}

\newcommand{\Pic}{\operatorname{Pic}}

\newcommand{\Sing}{\operatorname{Sing}}

\newcommand{\Bs}{\operatorname{Bs}}

\newcommand{\Cl}{\operatorname{Cl}}
\newcommand{\Aut}{\operatorname{Aut}}

\newcommand{\card}[1]{\operatorname{card}{(#1)}}

\newcommand{\rk}{\operatorname{rk}}
\newcommand{\rr}{\operatorname{r}}

\newcommand{\g}{\operatorname{g}}
\newcommand{\bb}{\operatorname{b}}

\newcommand{\Proj}{\operatorname{Proj}}

\newcommand{\CC}{\mathbb{C}}

\newcommand{\FF}{\mathbb{F}}
\newcommand{\QQ}{\mathbb{Q}}

\newcommand{\ZZ}{\mathbb{Z}}

\newcommand{\PP}{\mathbb{P}}

\newcommand{\OOO}{{\mathscr{O}}}

\newcommand{\NNN}{{\mathscr{N}}} 
\newcommand{\III}{{\mathcal{I}}}
\newcommand{\EEE}{{\mathscr{E}}}

\newcommand{\fU}{\mathfrak{U}}
\newcommand{\fY}{\mathfrak{Y}}

\newcommand{\fC}{\mathfrak{C}}
\newcommand{\fL}{\mathfrak{L}}
\newcommand{\fD}{\mathfrak{D}}
\newcommand{\fH}{\mathfrak{H}}
\newcommand{\fE}{\mathfrak{E}}
\newcommand{\fX}{\mathfrak{X}}

\newcommand{\f}{\mathfrak{f}}

\renewcommand{\emptyset}{\varnothing}

\newcommand{\type}[1]{$\mathbf{(#1)}$}

\newcommand{\qq}{\mathbin{\sim_{\scriptscriptstyle{\QQ}}}}

\newcommand{\comment}[1]{}

\newcommand{\xref}[1]{{\rm \ref{#1}}}

\newtheorem{theorem}{Theorem} 

\numberwithin{theorem}{section}
\numberwithin{equation}{theorem}

\newtheorem{mtheorem}[theorem]{}
\newtheorem{stheorem}[equation]{}

\theoremstyle{definition}
\newtheorem{case}[theorem]{}
\newtheorem{scase}[equation]{}

\newcounter{NN}
\renewcommand{\theNN}{{\rm\arabic{NN}${}^o$}}
\def\nr{\refstepcounter{NN}{\theNN}}

\newcounter{NO}
\renewcommand{\theNO}{{\rm(\Roman{NO})}}
\def\NR{\refstepcounter{NO}{\theNO}}

\begin{document}

\title{Singular Fano threefolds of genus 12}
\author{Yuri Prokhorov}
\thanks{
The author is supported by the RSF, grant No. 14-21-00053.
}

\address{
Steklov Mathematical Institute, Russia
\newline\indent
Department 
of Algebra, Moscow State
University, Russia
\newline\indent
National Research University Higher School of Economics, Russia
}
\email{prokhoro@mi.ras.ru}

\begin{abstract}
We study singular Fano threefolds of type $V_{22}$. 
\end{abstract}

\maketitle

\section{Introduction}
All Fano varieties in this paper are supposed to be three-dimensional, to have at worst 
terminal Gorenstein singularities and defined over an algebraically closed field of characteristic zero.
A Fano threefold $X$ is said to be \emph{of the main series} if the canonical class $K_X$
generates the Picard group $\Pic(X)$.
For a Fano threefold $X$ of the main series we can write $-K_X^3=2\g(X) -2$, where $\g(X)$ is an integer
which is called the \emph{genus} of $X$. It is known that $\g(X)$ takes the following values:
$\g(X)\in \{2,3,\dots,10,12\}$ (see \cite{Iskovskikh-Prokhorov-1999} and \cite{Namikawa-1997}).
Smooth Fano threefolds of the main series and genus $12$ were described by Iskovskikh \cite{Iskovskih1978a}
and Mukai \cite{Mukai-1989}. 
In this paper we study singular Fano threefolds of genus 12.
An important invariant of a Fano variety is $\rr(X):=\rk \Cl(X)$,
the rank of the Weil divisor class group.
The case $\rr(X)=1$ is already known:
\begin{mtheorem}{\bf Theorem \cite{Mukai-2002}, \cite{Prokhorov-planes}.}
Let $X$ be a $\QQ$-factorial Fano threefold of the main series with $\g(X)=12$. Then $X$ is smooth.
\end{mtheorem}
In \S \ref{section_Case_r=2} we prove the following result.
\begin{mtheorem}{\bf Theorem.}\label{Theorem-v22-rho-2}
Let $X$ be a Fano threefold of the main series with $\g(X)=12$.
Assume that $X$ satisfies one of the following equivalent conditions 
\begin{enumerate}
\item\label{Theorem-v22-rho-2-assumption-1}
the singular locus of $X$ consists of one ordinary double point,
\item\label{Theorem-v22-rho-2-assumption-2}
$\rr(X)=2$.
\end{enumerate}
Then $X$ is the midpoint of a Sarkisov link, i.e. it suits to the following commutative diagram 
\begin{equation} \label{(4.1.1)}
\vcenter{
\xymatrix{
&Y\ar[dl]_f\ar[dr]^{\pi}\ar@{-->}[rr]^{\chi}&&Y^+\ar[dl]_{\pi^+}\ar[dr]^{f^+}
\\
Z&&X&&Z^+
}
}
\end{equation} 
where $\pi$ and $\pi^+$
are small $\QQ$-factorializations, and $\chi$ is a flop. The morphisms 
$f$ and $f^+$ are extremal Mori contractions described as follows
\begin{center}
\begin{tabular}{l|l|p{150pt}|l|p{145pt}}
{\rm\tiny No.} & $Z$ & $f$ & $Z^+$ & $f^+$ 
\\[3pt]
\hline
&&&
\\[-4pt]
{\tiny\NR \label{theorem-v22-rho-2-P3}}&$\PP^3$&
the blowup of a smooth rational quintic curve 
$B\subset \PP^3$ which is not contained in a quadric
& $\PP^3$&the blowup of a smooth rational quintic curve which is not contained in a quadric
\\
{\tiny\NR\label{theorem-v22-rho-2-Q}} & $Q$ & the blowup of 
a smooth rational quintic curve $B\subset Q$ 
which is not contained in a hyperplane section& $\PP^2$ & a conic bundle whose
discriminant curve of degree $3$
\\
{\tiny\NR\label{theorem-v22-rho-2-V5}} & $V_5$ & the blowup of 
a smooth rational quartic curve $B\subset V_5 $ & $\PP^1$ & a del Pezzo fibration of degree $6$
\\
{\tiny\NR\label{theorem-v22-rho-2-DP5}} & $\PP^2$& $\PP_{\PP^2}(\EEE)\to \PP^2$, where $\EEE$ is a stable rank-$2$ vector bundle
on $\PP^2$ with $c_1=0$, $c_2=4$ 
& $\PP^1$ & a del Pezzo fibration of degree $5$
\end{tabular}
\end{center}
where $Q$ is a smooth quadric in $\PP^4$ and $V_5\subset \PP^6$ is a smooth del Pezzo threefold 
of degree $5$. 
\end{mtheorem}

Note that in general case the rank of the group $\Cl(X)$ for a Fano threefold of genus  12 
can be quite large:
in Example \ref{example-rho=8} we have $\rr(X)=8$ and moreover $\rr(X)=10$ for some toric varieties.
Proposition \ref{Proposition-le10} proves that $\rr(X)\le 10$  in the case where
$X$ does not contain planes. A sharp bound of the rank of the group $\Cl(X)$ is not known.

In principle, the proof of Theorem \ref{Theorem-v22-rho-2}
can be extracted from the series of papers \cite{Jahnke-Peternell-Radloff-II},
\cite{Takeuchi-2009}, \cite{Cutrone-Marshburn}, \cite{Blanc-Lamy-2012}.
However, all these papers consider the situation where $\QQ$-factorializations 
in \eqref{(4.1.1)} are smooth. 
We prefer to give a relatively short self-contained proof in general 
case.\footnote{The variety \ref{theorem-v22-rho-2-DP5} was erroneously omitted in 
\cite{Jahnke-Peternell-Radloff-II}}
Note that all the varieties in Theorem \ref{Theorem-v22-rho-2}
admit a moving decomposition \cite[\S 7]{Mukai-2002}.
Thus one cannot expect that Mukai's techniques \cite{Mukai-1989, Mukai-2002} works in our case. 

In \S\S \ref{section-I-II-III}-\ref{section-IV} we study cases \ref{theorem-v22-rho-2-P3}-\ref{theorem-v22-rho-2-DP5}
in details. In particular, we show that these cases occurs 
and describe general members of the families.
According to Mukai the moduli space $\mathcal M_{\mathrm{Fano}}^{12}$
of smooth Fano threefolds is $6$-dimensional and birational to the 
moduli space $\mathcal M_3$ of genus $3$ curves. 
Our four families \ref{theorem-v22-rho-2-P3}-\ref{theorem-v22-rho-2-DP5} are parametrized by 
varieties $\mathcal M_{\mathrm{Fano}}^{\ref{theorem-v22-rho-2-P3}}$,\dots, 
$\mathcal M_{\mathrm{Fano}}^{\ref{theorem-v22-rho-2-DP5}}$ 
of dimension $5$.

The original motivation of this work was to study $G$-varieties and finite subgroups of the Cremona group
(see \cite{Prokhorov2009e,Prokhorov-2013-im}). 
An algebraic variety is said to be 
a \textit{$G$-variety} if it is equipped with an action $G\to \Aut_{\Bbbk}(X)$ of a finite group $G$.
A projective $G$-variety $X$ is a 
\emph{$G$-Fano variety} if $X$ has at worst terminal 
singularities,
$-K_X$ is an ample Cartier divisor, and the rank of the 
invariant part $\Cl(X)^G$ of the Weil divisor class group 
equals $1$. 
$G$-Fano threefolds which are not of the main series were classified in \cite{Prokhorov-GFano-1,Prokhorov-GFano-2}.

The main result of this paper is the following.

\begin{mtheorem}{\bf Theorem.}\label{Theorem-main}
Let $X$ be a $G$-Fano threefold of the main series with $\g(X)=12$. Then $X$ either is smooth
or its singular locus consists of one ordinary double point.
In the latter case $X$ is the anticanonical image of 
the blowup of $\PP^3$ along a smooth rational curve of degree $5$
\textup(see Theorem \xref{Theorem-v22-rho-2} \xref{theorem-v22-rho-2-P3}\textup). 
\end{mtheorem}
In particular, the diagram \eqref{(4.1.1)} in the case \ref{theorem-v22-rho-2-P3}
gives a well known Cremona map $\PP^3\dashrightarrow \PP^3$ of degree $(3,3)$
(see \cite[Ch. VIII, Ex. 8, P. 185]{Semple-Roth-1949}).

One can modify the definition of $G$-varieties ($G$-Fano varieties etc.) for the arithmetic case:
the variety $X$ is defined over a non-closed field $\Bbbk$ and 
$G$ is the absolute Galois group 
acting on $\bar X=X\otimes \bar \Bbbk$ through the second factor.
Our result can be applied in this situation with small modifications.

Another motivation of for study singular Fano threefolds comes from 
some questions about affine and Moishezon non-projective varieties
(cf. \cite{Furushim2006}).
In \S \ref{section-I-II-III} we construct a new compactification of the 
affine space $\mathbb A^3$ as a variety of type \xref{theorem-v22-rho-2-V5}
(see Theorem \ref{compactification}). 

\par\bigskip
Part of the work was written during the author's stay 
at the Max-Planck- f\"ur Mathematik (Bonn).
The author would like to thank this institute for invitation, excellent 
working conditions and hospitality.
The author is also grateful to the referee for numerous comments which 
helped me to improve the manuscript.

\section{Preliminaries}
\textbf{Standard notation.}
\begin{itemize}
\item 
$\Cl(X)$ denotes the Weil divisor class group of a variety $X$,
\item 
$\rr(X)=\rk \Cl(X)$,
\item 
$\Pic(X)$ denotes the Picard group of a variety $X$,
\item 
$\uprho(X)$ is the Picard number of $X$,
\item 
$Q=Q_d\subset \PP^{d+1}$ is a smooth quadric of dimension $d$,

\item 
$Q'=Q'_d\subset \PP^{d+1}$ is a singular irreducible quadric of dimension $d$,
\item 
$V_5\subset \PP^6$ is a smooth quintic del Pezzo threefold,

\item $\NNN_{C/V}$ is the normal bundle of $C$ in $V$,

\item
$\FF_n=\PP_{\PP^1}(\OOO_{\PP^1}\oplus \OOO_{\PP^1}(n))$ is a rational ruled surface (Hirzebruch surface),

\item
$\iota(X)$ is the Fano index of a (generalized) Fano variety \cite[\S 2.1]{Iskovskikh-Prokhorov-1999}.

\item
$\langle B\rangle\subset \PP^n$ is the linear span of $B\subset \PP^n$.
\end{itemize}
A \emph{contraction} is a morphism with connected fibers of normal varieties. 
\par\smallskip\noindent
\textbf{Convention.}
Throughout this paper all (generalized) Fano varieties are supposed 
to have at worst terminal Gorenstein singularities. A
Fano threefold $X$ of Fano index $2$ is called a \emph{del Pezzo threefold}
(see e.g. \cite{Iskovskikh-Prokhorov-1999, Prokhorov-GFano-1}).

We will use systematically the following very useful observation.

\begin{mtheorem}{\bf Proposition \cite[Lemma 5.1]{Kawamata-1988-crep}.}
\label{Proposition-factorialization}
Let $U$ be a threefold with terminal Gorenstein singularities.
Then any $\QQ$-Cartier divisor on $U$ is Cartier.
\end{mtheorem}

\begin{mtheorem}{\bf Proposition \cite[Corollary 4.5]{Kawamata-1988-crep}.}
\label{Proposition-factorial}
Let $U$ be a threefold with terminal singularities.
Then there exists a small morphism $\pi: \tilde U\to U$ such that 
$\tilde U$ has only terminal $\QQ$-factorial singularities and 
$K_{\tilde U}=\pi^*K_U$.
\end{mtheorem}
If in the above notation the variety $U$ is Gorenstein, then 
$\tilde U$ has only terminal factorial singularities.
In this case, we say that $\pi$ is \emph{factorialization} of
$U$.

The classification of extremal contractions of terminal Gorenstein threefolds 
is almost the same as in the smooth case:
\begin{mtheorem}{\bf Theorem (\cite{Cutkosky-1988}).}
\label{theorem-Cutkosky}
Let $V$ be a threefold with terminal factorial singularities 
and let $f : V \to W$ be an extremal Mori contraction with $\dim W>0$. 
Then one of the following holds.

\setenumerate[0]{leftmargin=17pt,itemindent=12pt}
\begin{enumerate}
\item[\type{e_1}] 
$f$ is birational and contracts a surface $E$ to an \emph{irreducible} curve $B$. 
In this case $B$ has locally planar singularities and 
is contained in the smooth locus of $W$. The contraction $f$ is the blowup of the ideal of $B\subset W$
and the variety $W$ also has only terminal factorial singularities.

\item[\type{e_{2\text{-}5}}]
$f$ is birational and contracts a surface $E$ to a point $P$.
In this case $f$ 
is the blowup of the maximal ideal of $P\in W$ and one of the following holds:
{\rm
\begin{equation}\label{table-extremal-contractions}
\renewcommand{\arraystretch}{1.4}
\begin{array}{c|c|c|c|c|c|c}
&f(E) & E &\OOO_E(E) & \alpha_E&\updelta & K_V^2 \cdot E \\
\hline 
\text{\type{e_2}} & \text{smooth} & \PP^2 & \OOO_{} (-1) &2& 8 & 4
\\
\text{\type{e_{3\text{-}4}}} & \text{cA-point} &\text{$Q_2$ or $Q_2'$} & \OOO_{} (-1) &1& 2 & 2
\\ 
\text{\type{e_5}} & \frac12 (1,1,1) & \PP^2 & \OOO_{} (-2) &\frac12& \frac12 & 1
\end{array}
\end{equation}
}
\noindent
where $\updelta:=(-K_W)^3 - (-K_V)^3$
and $\alpha_E$ is the discrepancy of $E$.

\item[\type{d}]
$W$ is a smooth curve and $f$ is a del Pezzo fibration;
\item[\type{c}]
$W$ is a smooth surface and $f$ is a conic bundle.
\end{enumerate}
\end{mtheorem}
\begin{scase}
\label{remark-conic-bundles}
In the case \type{c} there exists a reduced curve $\Delta_f\subset W$
such that $f$ is a smooth morphism over $W\setminus \Delta_f$ and 
any fiber $V_w:=f^{-1}(w)$ over $w\in \Delta_f$ is a singular (either reducible or non-reduced) conic.
We say that $\Delta_f$ is the \emph{discriminant curve} of $f$.
If $\Delta_f=\emptyset$, then $f$ is a $\PP^1$-bundle, in particular, $V$ is smooth
(see e.g. \cite[Prop. 5.2]{Prokhorov-planes}).
If furthermore, the surface $W$ is rational, then this bundle is locally trivial.
\end{scase}

\begin{scase}\label{remark-DP-fibrations}
In the case \type {d} the general fiber $V_\nu$ is a smooth del Pezzo surface
of degree $d=K_{V_\nu}^2$ with $1\le d\le 9$, $d\neq 7$.
If $d=9$, then $f$ is a locally trivial $\PP^2$-bundle.
If $d=8$, then the general fiber $V_\nu$ is a smooth quadric.
\end{scase}

\begin{scase}
\label{length-extremal-contractions}
There exists the following natural exact sequence
\begin{equation*}
0 \longrightarrow \Pic(W) \overset{f^*}\longrightarrow \Pic(V) \overset{\cdot l}\longrightarrow \ZZ 
\end{equation*}
where $l$ is any curve contracted by $f$. To describe the cokernel of the right-side map, 
we denote 
\begin{equation}
\label{equation-length-extremal-contractions}
\mu_f:=\min \left\{ -K_V\cdot l \mid \text{$l\subset V$ is a curve contracted by $f$}\right\}.
\end{equation}
Thus the cokernel of the map $\Pic(V) \to \ZZ$ is a cyclic group of order $\mu_f$.
Then $\mu_f=1$ except for the following cases:
\begin{itemize}
\item
[]$\mu_f=2$: 
\type{e_2}, \type{c} with $\Delta_f=\emptyset$, \type{d} with 
$K_{V_\eta}^2=8$;
\item
[]$\mu_f=3$: \type{d} with $K_{V_\eta}^2=9$.
\end{itemize}
\end{scase}

\begin{mtheorem}{\bf Lemma \cite[Lemma 3.1]{Kaloghiros2011}, \cite[Proposition 5.1]{Prokhorov-planes}.}
\label{equation-E2-5-main}
Let $f: V\to W$ be the blowup of a reduced \textup(but possibly reducible\textup) 
locally planar curve $B$ as in 
\type{e_1} of Theorem \xref{theorem-Cutkosky}. Then 
\begin{equation}
\label{equation-E1-main}
\begin{array}{rcl}
(-K_V)^3 &=& (-K_W)^3 + 2K_W \cdot B + 2p_a(B)-2,
\\[6pt]
(-K_V)^2 \cdot E &=& -K_W\cdot B - 2p_a(B)+2,
\\[6pt]
(-K_V) \cdot E^2 &=& 2p_a(B)-2.
\end{array}
\end{equation}
Therefore,
\begin{equation}\label{equation-E1-main-0}
(-K_W)^3 - (-K_V)^3 = 2K_V^2 \cdot E + 2p_a(B)-2. 
\end{equation}
\end{mtheorem}

\begin{mtheorem}{\bf Lemma.}\label{Lemma-degree-increase}
Let $V$ be a projective threefold with at worst terminal $\QQ$-factorial singularities 
such that $-K_V$ is nef.
Let $f: V\to W$ be a birational contraction such that 
$-K_W$ is a nef $\QQ$-Cartier divisor.
Let $S\subset V$ be an irreducible surface which is not contained in the 
$f$-exceptional locus. Then 
\begin{equation*}
K_{V}^2\cdot S \le K_{W}^2\cdot f( S).
\end{equation*}
\end{mtheorem}
\begin{proof}
Write 
\begin{equation*}
K_V\qq  f^*K_W +\sum \alpha_i E_i, 
\end{equation*}
where $E_i$ are prime exceptional divisors, 
and $\alpha_i\in \QQ$. Since $-K_V$ is $f$-nef, $\alpha_i\ge 0$.
Put $E=\sum \alpha_i E_i$.
Then
\begin{multline*}\label{multline-induction-a}
K_{V}^2\cdot S = K_{V}\cdot(f^*K_{W}+E) \cdot S\le
K_{V}\cdot f^*K_{W} \cdot S=
\\
=
(f^*K_{W}+E)\cdot f^*K_{W}\cdot S 
\le (f^*K_{W})^2 \cdot S= K_{W}^2\cdot f( S). 
\end{multline*}
\end{proof}

\begin{case}{\bf Definition.}
A \emph{generalized Fano variety} is a projective variety $X$ with 
at worst terminal Gorenstein singularities such that the 
anticanonical divisor $-K_X$ is nef and big and the
natural morphism 
\begin{equation*}
\Phi: X \longrightarrow\bar X:=\Proj \oplus_{n\ge 0} H^0(X,-nK_X)
\end{equation*} 
to the anticanonical model
does not contract divisors. 
\end{case}
Note that in this situation $\bar X$
is a Fano variety (with 
at worst terminal Gorenstein singularities).
Conversely, if $\bar X$ is a Fano variety (as above) and 
$\Phi: X\to \bar X$ is its small factorialization (see Proposition \ref{Proposition-factorial}), then $X$ is a 
generalized Fano variety.

\begin{case}
Let $X$ be a generalized Fano threefold.
Put $\g(X):=-K_X^3/2+1$. 
By Riemann-Roch and Kawamata-Viehweg vanishing $\dim |-K_X|=\g(X)+1$.
Hence $\g(X)$ is an integer. It is called the \emph{genus} of $X$. 
The Picard group $\Pic(X)$ and the Weil divisor class group $\Cl(X)$ 
are finitely generated and torsion free \cite[Proposition 2.1.2]{Iskovskikh-Prokhorov-1999}. 
Moreover, there exists a natural embedding $\Pic(X) \hookrightarrow \Cl(X)$
as a primitive sublattice (see Proposition \ref{Proposition-factorialization}). 
The \emph{Fano index} of $X$ is the maximal integer $\iota=\iota(X)$ such that 
$K_X$ is divisible by $\iota$ in $\Pic(X)$.
The \emph{degree} of a surface $S$ in a generalized Fano threefold $X$ is the 
degree with respect to the anticanonical divisor, i.e. the
(positive) number $(-K_X)^2\cdot S$. We say that a surface 
$S\subset X$ is a \emph{plane} if its degree equals $1$. 
Note that for a Fano threefold $X$ of the main series
with $\g(X)\ge 4$,
the linear system $|-K_X|$ is base point free \cite{Jahnke-Radloff-2006}.
Hence, in this situation, any surface of degree $1$ is
isomorphic to $\PP^2$, i.e. it is a plane in the usual sense.
\end{case}

\begin{mtheorem}{\bf Lemma.}\label{lemma-vector-bundle-2}
Let $\EEE$ be a rank-$2$ vector bundle on a smooth surface $Z$,
let $Y:=\PP_{Z}(\EEE)$, and let 
$M$ be the tautological divisor on $Y$.
Then the following relations hold
\begin{gather}
\label{equation-vector-bundle-K}
-K_{Y}\sim 2M +f^*(-K_Z-c_1(\EEE)), 
\\
\label{equation-vector-bundle-Hirsch}
M^2=M\cdot f^* c_1(\EEE)-f^*c_2(\EEE),\quad M^3=c_1(\EEE)^2-c_2(\EEE). 
\\
\label{equation-vector-bundle-K3}
-K_{Y}^3=6K_Z^2+2c_1(\EEE)^2-8c_2(\EEE).
\end{gather}
\end{mtheorem}
\begin{proof}
Use the relative Euler exact sequence and the Hirsch formula.
\end{proof}

\section{MMP on generalized Fano threefolds}

Our proof of Theorem \ref{Theorem-main} is essentially uses the following result.
\begin{mtheorem}{\bf Theorem \cite{Prokhorov-planes}.}
\label{Theorem-planes}
Let $X$ be a $G$-Fano threefold of the main series with
$\g(X)\ge 6$. Then $X$ contains no any planes. 
\end{mtheorem}
It turns out that the absence of planes is important for 
application of the minimal model program (MMP) in our situation.
We explicitly describe steps of the MMP
on generalized Fano threefolds. This techniques was developed 
\cite{Prokhorov-2005a} and \cite{Kaloghiros2011}.

\begin{case}{\bf Set-up.}\label{setup}
Let $X$ be a Fano threefold (with at worst terminal Gorenstein singularities).
Assume that $X$ does not contain planes 
and $\rr(X)>1$. 
Let $\pi: \tilde X\to X$ be a 
factorialization (see Proposition \ref{Proposition-factorial}). 
Then $K_{\tilde X}=\pi^*K_X$ and so $\tilde X$ is a generalized Fano threefold with $\uprho(\tilde X)=\rr(X)>1$.
\end{case}

The following important lemma shows that the class of generalized $\QQ$-factorial Fano threefold
containing no planes is closed under the MMP. 

\begin{stheorem}{\bf Lemma \cite{Prokhorov-2005a}, \cite{Kaloghiros2011}.}
\label{lemma-claim-mmp}
Let $V$ be a generalized $\QQ$-factorial Fano threefold.
Assume that $V$ does not contain planes.
Let $\varphi: V \dashrightarrow W$ be either a divisorial Mori contraction or a flop.
Then $W$ is also a generalized $\QQ$-factorial Fano threefold
and $W$ does not contain planes.
\end{stheorem}

\begin{proof}
The assertion is obvious if $f$ is a flop, so we assume that $f$ is a divisorial Mori contraction.
Since $V$ contains no planes, the contraction $f$ is not of type \type{e_5}.
Then by Theorem \ref{theorem-Cutkosky} the variety $W$ has only factorial terminal singularities.
Moreover, $-K_W$ is nef and big by \cite[Prop. 4.5]{Prokhorov-2005a}
(note that $f$ 
cannot be a ``bad'' contraction of type $(2,0)_o^-$, $(2,1)_{o0}^-$,$(2,1)_{o1}^-$
because the exceptional divisor is not a plane).
Finally, $W$ contains no planes by Lemma \ref{Lemma-degree-increase}.
\end{proof}

\begin{case}\label{MMP}
Let $\Theta$ be any divisor on $\tilde X$. According to 
\cite[2.6-2.8]{Prokhorov-Shokurov-2009} we can run the $\Theta$-MMP:
\begin{equation*}
\varphi: \tilde X=\tilde X_0\dashrightarrow 
\cdots
\dashrightarrow
\tilde X_N.
\end{equation*}
By Lemma \ref{lemma-claim-mmp} 
on each step $\varphi_k: \tilde X_k \dashrightarrow \tilde X_{k+1}$
the following assertions hold.
\begin{enumerate}
\item \label{Claim-1}
$\tilde X_k$ has only factorial terminal singularities; 
\item \label{Claim-2}
$\tilde X_k$ is a generalized Fano threefold;
\item \label{Claim-3}
$\tilde X_k$ does not contain planes;
\item\label{Claim-4}
$\varphi_k$ is either a flop or an extremal Mori contraction of type \type{e_1}-\type{e_4}.
\end{enumerate}
Note that in the ``classical'' case $\Theta=K_{\tilde X}$ all the 
steps $\varphi_k$ must be divisorial contractions of type 
\type{e_1}-\type{e_4}, in particular, $\varphi$ is a morphism.

\begin{scase}\label{sMMP}
We get the following diagram:
\begin{equation}
\label {equation-MMP}
\vcenter{
\xymatrix@C+15pt{
&\llap{$\varphi:$ }\tilde X=\tilde X_0\ar@{-->}[r]^{\varphi_1}\ar@<-2.8ex>[d]^\pi&\tilde X_1\ar[d]^{\pi_1}\ar@{-->}[r]^{\varphi_2}
&\dots\ar@{-->}[r]^{\varphi_N} &\tilde X_{N}\ar[d]^{\pi_N}&
\\
&X=X_0&X_1&&X_{N}&
} }
\end{equation}
where $X_k$ is the pluri-anticanonical
image of $\tilde X_k$. Denote by $\Theta^{(k)}$
the proper transform of $\Theta$ on $\tilde X_k$. 
\footnote{\label{footnote-Kaloghiros2011} 
\cite{Kaloghiros2011} asserts that $\uprho(X_i)=1$
for all $i$.
This is wrong in general.
Indeed, let $X$ be a quintic del Pezzo threefold with $\rr(X)=4$ \cite[\S 7]{Prokhorov-GFano-1}.
For a suitable choice of factorialization, we have $N=3$, $\tilde X_3=X_3=\PP^3$,
and $\varphi_i$ are contractions of type \type{e_2}.
Then $X_1$ is a sextic del Pezzo threefold with $\rr(X)=3$ and $\uprho(X_1)=2$.
}
At the end we get one of the following possibilities:
\begin{itemize}
\item 
There exists a $\Theta^{(N)}$-negative Mori contraction
$\upsilon: \tilde X_N\to Z$  to a variety $Z$
of dimension $<3$. In particular,  
$\uprho(\tilde X_N/Z)=1$ and $-K_{\tilde X_N}$ is $\upsilon$-ample.
If $\dim Z>0$, then $\upsilon$ is a  contraction of type \type{c} or \type{d}.
\item
The divisor $\Theta^{(N)}$ is nef.
\end{itemize}
\end{scase}
\end{case}

\begin{scase}{\bf Definition.}
A \emph{weak del Pezzo surface} is a projective surface with at worst Du Val singularities such that
the anticanonical divisor $-K_X$ is nef and big. 
\end{scase}
\begin{scase}{\bf Remark.}\label{Remark-weak-del-Pezzo}
It is easy to see that the class of weak del Pezzo surfaces is closed under 
birational contractions.
\end{scase}

\begin{mtheorem}{\bf Lemma.}\label{MMp-to-surface}
Let $X$ be a generalized Fano threefold
and let $\upsilon: X\to Y$ be a contraction to a surface. 
Assume that $X$ does not contain planes.
Then the following 
assertions hold:
\begin{enumerate}
\item \label{MMp-to-surface-i}
$Y$ is a weak del Pezzo surface.
\item \label{MMp-to-surface-ii}
If $X$ is $\QQ$-factorial and $\upsilon$ is an extremal Mori contraction, then 
$Y$ is a smooth del Pezzo surface.
\item \label{MMp-to-surface-iii}
If under the assumptions of \ref{MMp-to-surface-ii} the surface 
$Y$ contains a $(-1)$-curve $\Gamma$, then 
there exists the following diagram:
\begin{equation*}
\xymatrix@R=15pt{
X\ar@{-->}[r]^{\chi}\ar[d]^{\upsilon}&X^+\ar[r]^{\varphi}&X'\ar[d]^{\upsilon'}
\\
Y\ar[rr]&&Z
} 
\end{equation*}
where $Z$ is a smooth surface, $Y\to Z$ is the contraction of $\Gamma$, $\chi$ is either an isomorphism or a flop, 
$\varphi$ is an extremal divisorial contraction, and 
$\upsilon'$ is an extremal contraction of type \type{c}. 
\end{enumerate}
\end{mtheorem}

\begin{proof}
\ref{MMp-to-surface-i}
As in \eqref{equation-MMP}, run the $K$-MMP on $X$ over $Y$:
\begin{equation}\label{equation-diggram-mmp-conic-bundle}
\vcenter{
\xymatrix@R=15pt@C=50pt{
X\ar[d]^{\upsilon}&\tilde X\ar[l]_{\pi}\ar[r]^{\varphi}&\hat X\ar[d]^{\tau}& 
\\
Y&&\hat Y\ar[ll]
} 
}
\end{equation}

Here $\pi$ is a suitable factorialization,
$\varphi$ is a composition of divisorial contractions, and 
$\tau$ is an extremal contraction of type \type{c}.
By Lemma \ref{lemma-claim-mmp} $\hat X$ is a 
generalized Fano threefold.
Then by \cite[Proposition 5.2 (i)]{Prokhorov-2005a}
$\hat Y$ is a \emph{smooth} weak del Pezzo surface. 
This implies that $Y$ is a weak del Pezzo surface (see Remark \ref{Remark-weak-del-Pezzo}).

\ref{MMp-to-surface-ii} and \ref{MMp-to-surface-iii}
Again by Theorem \ref{theorem-Cutkosky}\type{c} and \cite[Proposition 5.2 (i)]{Prokhorov-2005a} 
$Y$ is a smooth weak del Pezzo surface.
Let $\Gamma\subset Y$ be a curve with $\Gamma^2<0$ and let $\tau: Y\to U$ be its contraction.
Clearly, $\uprho(X/U)=2$ and so the relative Mori cone $\overline{\operatorname{NE}}(X/U)$
has two extremal rays, say $\mathcal R_1$ and $\mathcal R_2$.
Put $F:=\upsilon^*(\Gamma)$.
We may assume that $\mathcal R_1$ is generated by the fibers of $\upsilon$
and so $F\cdot \mathcal R_1=0$. Since $\Gamma^2<0$, 
for any curve $\Gamma'\subset X$ dominating $\Gamma$ we have 
$F\cdot \Gamma'<0$. Hence 
the divisor $F:=\upsilon^*(\Gamma)$ is not nef and $F\cdot \mathcal R_2<0$. 
Thus we can run the $F$-MMP over $U$ starting from $\mathcal R_2$:
\begin{equation*}
\xymatrix@R=15pt@C=50pt{
X\ar@{-->}[r]^{\chi}\ar[d]^{\upsilon}&X^+\ar[r]^{\varphi}&X'\ar[d]^{\upsilon'}
\\
Y\ar[r]^{\tau}&U&Y'\ar[l]_{\tau'}
} 
\end{equation*}
Here $\chi$ is either an isomorphism or a flop in $\mathcal R_2$ and
$\varphi$ is either an isomorphism or a divisorial contraction. 
By Lemma \ref{lemma-claim-mmp} $X'$ is a 
generalized Fano threefold.
Thus $\upsilon'$ is an extremal contraction to a 
smooth surface $Y'$ (see Theorem \ref{theorem-Cutkosky}\type{c}). 
Note that $\uprho(X^+)=\uprho(X)=\uprho(U)+2$.

Assume that $U$ is singular. Then $\uprho(Y')>\uprho(U)$ and so 
$\uprho(X^+)=\uprho(X')=\uprho(Y')+1$, i.e. $\varphi$ is an isomorphism.
Moreover, $Y'$ is a minimal resolution of $U$ and so $Y\simeq Y'$.
Since both $-K_{X}$ and $-K_{X'}$ are ample over $Y=Y'$, the map 
$\varphi \comp \chi$ must be an isomorphism, a contradiction.

Thus $\Gamma$ is a $(-1)$-curve and $U$ is smooth. 
In particular, this implies that $Y$ is a del Pezzo surface.
If $\varphi$ is an isomorphism, then $\tau'$ is 
the blowup of the (smooth) point $\tau(\Gamma)\in U$ and $Y'\simeq Y$. We get a contradiction as above.
Hence $\varphi$ is a (single) divisorial contraction.
%
\end{proof}

\begin{mtheorem}{\bf Lemma.}\label{Lemma-surface-rho=2}
One can run the MMP \eqref{equation-MMP} so that 
$\uprho(Z)\le 2$.
\end{mtheorem}

\begin{proof}
Assume that $\uprho(Z)\ge 3$. 
Then $Z$ must be a smooth rational surface (see Lemma \ref{MMp-to-surface}\ref{MMp-to-surface-ii}) and 
it contains a $(-1)$-curve. Then by Lemma \ref{MMp-to-surface}\ref{MMp-to-surface-iii}
we can run the MMP further
until we get a surface $Z_M$ with $\uprho(Z_M)\le 2$.
\end{proof}

\begin{mtheorem}{\bf Lemma.}\label{lemma-contractions}
Let $X$ be a generalized Fano threefold
and let $f: X\to Y$ be a birational contraction.
Assume that $X$ does not contain planes.
Then $Y$ is also a generalized Fano threefold.
\end{mtheorem}

\begin{proof} 
As in \eqref{equation-MMP}, run the $K$-MMP on $X$ over $Y$:
we start with a suitable factorialization $\pi: \tilde X\to X$
and after a number of divisorial contraction we obtain 
a minimal model $\tilde Y$ over $Y$.
By 
Lemma \ref{lemma-claim-mmp} $\tilde Y$ is a 
generalized Fano threefold.
Hence the morphism $\tilde Y\to Y$ is small and crepant.
So, $Y$ is a generalized Fano threefold as well.
\end{proof}

\begin{stheorem}{\bf Corollary.}\label{corollary-degree-difference}
In the notation of Lemma \xref{lemma-contractions} we have 
\begin{equation*}
\rr(X)-\rr(Y) \le \frac 12 K_X^2\cdot E,
\end{equation*}
where $E$ is the $f$-exceptional divisor.
\end{stheorem}

\begin{proof}
Since $X$ contains no planes, $E$ has at most $(K_{X}^2 \cdot E)/2$ components.
\end{proof}

\section{Deformations of Fano threefolds}

\begin{mtheorem}{\bf Theorem \cite{Namikawa-1997}.}\label{mtheoremNamikawa-1997}
Let $X$ be a generalized Fano threefold.
Then $X$ is smoothable, that is, there exists a flat family 
$\fX \to (\fU\ni 0)$ over a small disc
$(\fU\ni 0) \subset \CC$ such that $\fX_0\simeq 
X$ and a general fiber $\fX_u$, $u \in \fU\setminus \{0\}$, is a smooth generalized Fano
threefold. 
\end{mtheorem}

\begin{mtheorem}{\bf Theorem \cite{Jahnke2011}.}\label{mtheoremJahnke}
Let $X$ be a generalized Fano threefold
and let $\fX \to (\fU\ni 0)$ be its smoothing as above.
Then $\fX$ is normal and has at worst isolated terminal factorial singularities.
Moreover, there are natural identifications 
\begin{equation*}
\Pic(X) =\Pic(\fX_u ) = \Pic(\fX)
\end{equation*}
so that $K_{ \fX_ s} = K_X$.
\end{mtheorem}

\begin{scase}{\bf Remark.}\label{remark-plane-def}
If in the above notation for $u\neq 0$ the fiber $\fX_u$ contains a plane, then 
the same holds for the central fiber.
\end{scase}

\begin{stheorem}{\bf Corollary.}\label{Corollary-deformations-c}
Let $X$ be a Fano threefold with terminal Gorenstein singularities
such that $\uprho(X)=1$ and $-K_X^3>22$.
Then $\iota(X)>1$. If furthermore $-K_X^3\ge 40$, then $X$ is either $\PP^3$ or 
a quadric in $\PP^4$, or a quintic del Pezzo threefold in $\PP^6$. 
If $-K_X^3\ge 40$ and $X$ is $\QQ$-factorial, then $X$ is in fact smooth.
\end{stheorem}

\begin{proof}
Let $\fX \to \fU\ni 0$ be a smoothing as above.
A general fiber $\fX_u$ is a smooth Fano threefold with
$\uprho(\fX_u)=\uprho(X)=1$ and $-K_{\fX_u}^3=-K_X^3 >22$.
By the classification of smooth Fano threefolds with $\uprho=1$ we have $\iota(\fX_u)>1$
(see \cite[\S 12.2]{Iskovskikh-Prokhorov-1999})
and by Theorem \ref{mtheoremJahnke} $\iota(X)=\iota(\fX_u)>1$.
The case $\iota(X)\ge 3$ is well-known (see e.g. \cite[\S 3.1]{Iskovskikh-Prokhorov-1999}).
For the case $\iota(X)=2$ we refer to 
\cite[Corollary 8.7]{Prokhorov-GFano-1}.
\end{proof}

\begin{case}{\bf Notation.}\label{notation-deformations}
Below in this section we assume that 
$X$ is a Fano threefold with terminal Gorenstein singularities
and $X$ does not contain planes. 
Let $\fX \to (\fU\ni 0)$ be its $1$-parameter smoothing as above.
\end{case}

\begin{mtheorem}{\bf Proposition (cf. \cite[Proposition 6.3]{Prokhorov-GFano-2}).}
\label{proposition-nef-ample-deformations}
Notation as in \xref{notation-deformations}. Let $\fL$ be a divisor on $\fX$.
Then $\fL$ is nef \textup(resp. ample\textup) if and only if the restriction $\fL|_{\fX_u}$ 
is nef \textup(resp. ample\textup) for some $u\in \fU$.
\end{mtheorem}

\begin{mtheorem}{\bf Proposition (cf. \cite[Corollary 6.4]{Prokhorov-GFano-2}).} 
\label{Corollary-deformations}
Notation as in \xref{notation-deformations}.
Let $\f_u : \fX_u \to \fY_u$ be an extremal contraction. 
Then there exists an extremal
contraction $\f : \fX \to \fY$ over $\fU$ 
such that the restriction $\f|_{\fX_u}$ coincides with $\f_u$, 
where the variety $\fY$ 
is $\QQ$-factorial.
Let $f: X=\fX_0\to Y=\fY_0$ be the restriction 
of $\f$ to $X=\fX_0$. Then $Y$ is normal, $\QQ$-Gorenstein 
and $f$ has connected fibers.
Moreover, the following assertions hold:
\begin{enumerate}
\item \label{def-surfaces}
If $\dim \fY=3$, then $Y$ is a weak del Pezzo surface.
\item \label{def-birational}
If $\f$ is birational, then
$Y$ 
is a generalized Fano threefold.
If moreover, $\fY_u$ is smooth, then there are natural identifications 
$\Pic(Y) =\Pic(\fY_u )$ so
that $K_{ \fY_ s} = K_Y$, in particular, $\uprho(Y)=\uprho(\fY_ s)$ and $\iota(Y)=\iota(\fY_ s)$.
\end{enumerate}
\end{mtheorem}

\begin{proof}
The existence of $\f$ immediately follows from Proposition \ref{proposition-nef-ample-deformations}
(see \cite[Corollary 6.4]{Prokhorov-GFano-2}). 
By our assumptions $\fX$ contains no planes. 
Hence by \cite[Theorem 1.1]{Kachi1998} there are no flipping contractions on $\fX$.
This implies that the variety $\fY$ 
is $\QQ$-factorial (see \cite[Lemma 5-1-5, Proposition 5.1.6]{KMM})
and so $Y=\fY_0$ is $\QQ$-Gorenstein.
By the projection formula and Kawamata-Viehweg vanishing theorem
\begin{equation*}
R^1 \f_* \OOO_{\fX} (-X)= R^1 \f_* \OOO_{\fX}\otimes \f^*\OOO_{\fY}(-Y)=
R^1 \f_* \OOO_{\fX}\otimes \OOO_{\fY}(-Y)=0.
\end{equation*}
Applying $\f_*$ to the exact sequence
\begin{equation*}
0\longrightarrow \OOO_{\fX} (-X)\longrightarrow 
\OOO_{\fX} \longrightarrow \OOO_X \longrightarrow 0
\end{equation*}
we obtain
\begin{equation*}
0\longrightarrow \OOO_{\fY} (-Y)\longrightarrow 
\OOO_{\fY} \longrightarrow \f_* \OOO_X \longrightarrow 0.
\end{equation*}
Therefore, $\f_* \OOO_X=\OOO_Y$, $Y$ is normal, 
and $f$ has connected fibers, i.e. $f$ is a contraction.
Then apply Lemmas \ref{MMp-to-surface},
\ref{lemma-contractions} and Theorem \ref{mtheoremJahnke}.
\end{proof}

\begin{mtheorem}{\bf Lemma.}\label{Lemma-p1p2}
Notation as in \xref{notation-deformations}. 
Let $\f: \fX\to \fY$ be an extremal contraction of relative dimension one.
Assume that a general fiber $\fX_u$ has a contraction $\f_u': \fX_u\to \fY_u'\simeq \PP^1$
which is not passed through $\fY_u$.
Then $\f$ has no two-dimensional fibers.
\end{mtheorem}

\begin{proof}
Assume that $\f_0: X=\fX_0\to \fY_0$ has a two-dimensional fiber $F\subset X$.
By Proposition \ref{proposition-nef-ample-deformations} there exists a 
contraction $\f': \fX\to \fY'$ that extends 
$\f_u'$. Clearly, $\fY'\to \fU$ is a $\PP^1$-bundle.
Take ample divisors $\fH$ and $\fH'$ on $\fY$ and $\fY'$, respectively.
Let $\fL:=\f^* \fH$ and $\fL':=\f^* \fH'$.
Then the divisor $\fL_u+\fL'_u$ is ample by our assumption.
By Proposition \ref{proposition-nef-ample-deformations} the same holds for $\fL+\fL'$.
On the other hand, $\fL|_F=0$ and $\dim \f'(F)=1$.
Hence, the linear system $|\fL+\fL'|$ (i.e. the corresponding morphism) contracts $F$, a contradiction.
\end{proof}

\begin{mtheorem}{\bf Lemma.}\label{Lemma-def-surface-smooth}
Notation as in \xref{notation-deformations}.
Let $\f: \fX\to \fY$ be a contraction \textup(over 
$\fU\ni 0$\textup) all whose fibers are one-dimensional.
Then $\fY$ is smooth, $\f$ is a conic bundle, and its
discriminant locus $\fC\subset \fY$ is either empty or 
flat of relative dimension one over $\fU$.
Moreover, the fiber $\fY_0$ is a smooth weak del Pezzo surface of degree $K_{\fY_0}^2=K_{\fY_u}^2$
and the discriminant locus of $\f_0$ coincides with $\fC_0\subset \fY_0$.
If $\fC\subset \fY$ is empty, then $X$ is smooth and $\f_0$ is a $\PP^1$-bundle.
\end{mtheorem}

\begin{proof}
Similar to \cite[Theorem 7]{Cutkosky-1988} one can show that $\fY$ is smooth and $\f$ is a conic bundle. 
By Lemma \ref{MMp-to-surface} \ $\fY_0$ is a smooth weak del Pezzo surface.
The rest is obvious.
\end{proof}

\begin{mtheorem}{\bf Lemma.}\label{Lemma-product}
Notation as in \xref{notation-deformations}.
Assume that a general fiber $\fX_u$ 
is isomorphic to a product $\fX_u\simeq Y_u\times \PP^1$. 
Then the special fiber $X$ is smooth.
\end{mtheorem}

\begin{proof}
Consider the extremal contraction $\f: \fX\to \fY$ that corresponds to the projection
$\f_u:\fX_u\simeq Y_u\times \PP^1\to Y_u$. 
By Lemma \ref{Lemma-p1p2} the contraction $f: X\to Y$ has no two-dimensional fibers and
by Lemma \ref{Lemma-def-surface-smooth} 
the variety $X$ is smooth.
\end{proof}

\begin{mtheorem}{\bf Lemma.}\label{Lemma-Y12}
Notation as in \xref{notation-deformations}.
Assume that a general fiber $\fX_u$ 
is isomorphic to a divisor of bidegree $(1,2)$ in $\PP^2\times \PP^2$.
Then the special fiber $X$ has the same form
and $\rr(X)\le 5$.
\end{mtheorem}

\begin{proof}
Two projections $p_{u,i} : \fX_u\to \PP^2$
induce two contractions $p_i: X\to \PP^2$.
The divisors $H_i:= p_i^*\OOO_{\PP^2}(1)$ generate the Picard group $\Pic(X)$.
Moreover,
\begin{equation*}
-K_X=2H_1+H_2,\quad H_1^3=H_2^3=0,\quad H_1^2\cdot H_2=2,\quad H_1\cdot H_2^2=1.
\end{equation*}
Since $\uprho(X)=2$, the product map 
$\pi=p_1\times p_2: X\to \PP^2\times \PP^2$ is finite.
It is easy to see that $\pi$ is birational and the image $Y:=\pi(X)$ is a divisor 
of bidegree $(1,2)$. Then by the adjunction formula
$K_X=\pi^* K_Y$. Hence the map $\pi$ is an isomorphism onto its image.
If the projection $p_2: X\to \PP^2$ has a two-dimensional fiber, say $F$, then the anticanonical image of $F$ 
must be a plane. This contradicts our assumptions.
Hence $p_2$ is an equidimensional conic bundle with discriminant curve $\Delta_f\subset \PP^2$ of degree $3$.
Then $\Delta_f$ has at most $3$ components
and so $\rr(X)\le \rr(\PP^2)+1+ 3=5$.
\end{proof}

\begin{mtheorem}{\bf Lemma.}\label{Proposition-XN}
Let $X$ be a Fano threefold.
Assume that $X$ does not contain planes, has no 
any birational contractions, and $-K_X^3> 30$.
Then either $\iota(X)\ge 2$ or $X\simeq \PP^1\times \PP^2$.
\end{mtheorem}

\begin{proof}
Assume that $\iota(X)=1$.
Let $\fX \to
\fU\ni 0$ be a smoothing as in Theorem \ref{mtheoremNamikawa-1997}. 
Then for $0\neq u\in \fU$, 
the fiber $\fX_u$ is a smooth Fano threefold with $\uprho(\fX_u)=\uprho(X)$,
$\iota(\fX_u)=1$, and it has no 
any birational contractions. 
Then by \cite{Mori1981-82} we have only one possibility:
$\fX_u\simeq \PP^1\times \PP^2$.
By Lemma \ref{Lemma-product} \ 
$X\simeq \PP^1\times \PP^2$.
\end{proof}

\begin{mtheorem}{\bf Proposition.}\label{proposition-24}
Let $X$ be a Fano threefold with at worst terminal Gorenstein singularities 
and $-K_X^3=24$. Assume that $X$ contains no planes.
Then $\rr(X)\le 9$.
\end{mtheorem}

\begin{proof}
If $\uprho(X)=1$, then $\iota(X)>1$ and $X$ is a del Pezzo threefold
(see Corollary \ref{Corollary-deformations-c}).
In this case, 
$\rr(X)\le 6$ by \cite[Corollary 3.13]{Prokhorov-GFano-1}.
Thus we assume that $\uprho(X)>1$.
Consider a $1$-parameter smoothing $\fX \to (\fU\ni 0)$
as in Theorem \ref{mtheoremNamikawa-1997}.
A general fiber $\fX_u$ is a smooth Fano threefold with
$\uprho(\fX_u)=\uprho(X)>1$ and $-K_{\fX_u}^3=-K_X^3=24$.
Now we apply the classification \cite{Mori1981-82} to $\fX_u$.

If $\uprho(\fX_u)>4$, then $\fX_u$ 
is a product $\fX_u\simeq Y_u\times \PP^1$, where $Y_u$ is a del Pezzo surface of 
degree $4$ (see \cite[Table 5]{Mori1981-82}). 
Then by Lemma \ref{Lemma-product} the variety $X$ has the same form 
and so $\rr(X)=7$. Thus we may assume that $2\le \uprho(\fX_u)\le 4$.

If $\fX_u$ is not isomorphic to a blowup of a smooth Fano threefold, then 
$\fX_u$ is a double cover of $\PP^1\times\PP^2$
whose branch locus is a divisor of bidegree $(2,2)$ 
(\cite[Table 2, $\mathrm{n^o}$ $\mathrm{18^{o}}$]{Mori1981-82}).
In this case $\fX_u$ has two extremal contractions:
a conic bundle $\f_u: \fX_u\to \PP^2$ with discriminant curve $D_u\subset \PP^2$ of degree $4$
and a del Pezzo fibration $\f'_u: \fX_u\to \PP^1$.
Consider the extremal contraction $\f: \fX\to \fY$ that extends 
$\f_u$.
Then $Y=\fY_0$ is a del Pezzo surface with at worst Du Val singularities and $\uprho(Y)=1$.
Since $K_Y^2=K_{\fY_u}^2=9$, we have $Y\simeq \PP^2$.
By Lemma \ref{Lemma-p1p2} the contraction $f: X\to Y$ has no two-dimensional fibers
and $\f$ is a conic bundle.
Let $\fD\subset \fY$ be the discriminant curve of $\fX$.
Thus $D_u= \fY_u\cap \fD$. Moreover, $D_0= \fY_0\cap \fD$ is the discriminant curve of $\fX_0\to \fY_0\simeq \PP^2$
with $\deg D_0=\deg D_u=4$.
Then as in the proof of Lemma \ref{Lemma-Y12} we have
$\rr(X)\le \rr(Y)+1+ 4\le 6$.

In the remaining cases $\fX_u$ is isomorphic to a blowup of a smooth Fano threefold $\fY_u$
along a smooth curve $B_u$. 
Consider the extremal contraction $\f: \fX\to \fY$ that extends
$\fX_u\to \fY_u$. Let $Y:=\fY_0$.
By Proposition \ref{Corollary-deformations} 
$Y$ is a generalized Fano threefold with 
$\uprho(Y)=\uprho(\fY_u)$, $-K_Y^3=-K^3_{\fY_u}$, and $\iota(Y)=\iota(\fY_u)$.
Let $\fE\subset \fX$ be the exceptional divisor,
let $E_u:=\fE\cap \fX_u$, and let $E:=\fE_0$.
Obviously, the morphism $\fE\to \fU$ is flat and so
$K_X^2 \cdot E= K_{\fX_u}^2 \cdot E_u$.
By Corollary \ref{corollary-degree-difference},
\begin{equation}\label{equation-24}
\rr(X)\le \rr(Y)+ \frac12 K_{\fX_u}^2 \cdot E_u. 
\end{equation}
By \cite{Mori1981-82} 
there are the following possibilities: 
\begin{center}
\begin{tabular}{c|c|p{160pt}|c|c|c}\setcounter{NN}{0}
&$\uprho(\fX_u)$&$\fY_u$& $\rr(Y)$ &$p_a(B_u)$ &$K_{\fX_u}^2 \cdot E_u$\\
\hline
&&&
\\[-8pt]
\nr\label{nnext-rho=2a}
&$2$& a quadric $Q\subset\PP^4$ & $\le 2$ &$1$ &$15$
\\
\nr\label{nnext-rho=3a}
&$3$&a del Pezzo $3$-fold $V_6\subset \PP^7$& $\le 3$& $1$ &$12$
\\
\nr\label{nnext-rho=3b}
&$3$&$Y_{2,1}\subset \PP^2\times \PP^2$& $\le 5$ &$0$ & $4$
\\
\nr\label{nnext-rho=4}
&$4$&$\PP^1\times\PP^1\times\PP^1$& $3$
&$1$ &$ 12$\\
\hline
\end{tabular}
\end{center}
where $Y_{2,1}\subset \PP^2\times \PP^2$ is a divisor of bidegree $(2,1)$.
In the above table the value of $K_{\fX_u}^2 \cdot E_u$ is computed by using
\eqref{equation-E1-main-0}. To estimate $\rr(Y)$, we note that 
$Y$ is a (possibly singular) quadric in the case \ref{nnext-rho=2a},
$\rr(Y)\le 3$ in cases \ref{nnext-rho=3a} and \ref{nnext-rho=4}
by \cite[Corollary 3.13]{Prokhorov-GFano-1}, and $\rr(Y)\le 5$ in the case \ref{nnext-rho=3b}
by Lemma \ref{Lemma-Y12}.
Combining the above table with \eqref{equation-24} we get the desired estimate $\rr(X)\le 9$.
\end{proof}

\section{Proof of Theorem \xref{Theorem-v22-rho-2}.}
\label{section_Case_r=2}
In this section we prove Theorem \xref{Theorem-v22-rho-2}.

\begin{mtheorem}{\bf Proposition.}
Let $U$ be a Fano threefold whose singular locus consists of 
one ordinary double point. Then $\rr(U)\le \uprho(U)+1$. 
\end{mtheorem}
\begin{proof}
We may assume that $U$ is not $\QQ$-factorial.
Let $\pi: \tilde U\to U$ be a small factorialization (see Proposition \ref{Proposition-factorial}).
Then $\rr(\tilde U)=\rr(U)$. On the other hand,
$\uprho(\tilde U/U)=1$ because $\pi$ contracts exactly one
irreducible curve.
\end{proof}

Thus in Theorem \xref{Theorem-v22-rho-2} the implication 
\ref{Theorem-v22-rho-2-assumption-1}
$\Longrightarrow$ \ref{Theorem-v22-rho-2-assumption-2} holds and
we may assume that $\rr(X)=2$.

Let $X$ be a Fano threefold with 
$\g(X)=12$, $\uprho(X)=1$, and $\rr(X)=2$. 
Let $\pi:Y\to X$ be a small factorialization.
By our assumption $\pi$ is not an isomorphism and $Y$ is a \emph{generalized} Fano threefold 
with $\uprho(Y)=2$.
The existence of diagram \eqref{(4.1.1)} is a standard fact 
(see e.g. \cite[\S 4.1]{Iskovskikh-Prokhorov-1999}).
Consider the possibilities for the contraction $f$ according to Theorem \ref{theorem-Cutkosky}.
Note that our diagram \eqref{(4.1.1)} is symmetric.
So we can interchange $f$ and $f^+$.
In particular, if one of contractions $f$ and $f^+$ is birational,
then we may assume that this holds for $f$.

\begin{mtheorem}{\bf Lemma.}\label{lemma-r=2-e2}
Neither $f$ nor $f^+$ is of type \type{e_{2}}-\type{e_{3\text{-}4}}.
\end{mtheorem}
\begin{proof}
Assume that $f$ is of type \type{e_2} or \type{e_{3\text{-}4}}.
Let $E\subset Y$ be the exceptional divisor.
Then $Z$ is a $\QQ$-factorial Fano threefold with $\uprho(Z)=1$.
By \eqref{table-extremal-contractions} we have $-K_Z^3=30$ or $24$.
Hence, $\iota(Z)\ge 2$. 
Note that $-K_{Z}^3$ must be divisible by $\iota(Z)^3$.
The only possibility is $\iota(Z)= 2$ and $-K_Z^3=24$,
i.e. $Z$ is a del Pezzo threefold of degree $3$ and $f$ is of type \type{e_{3\text{-}4}}.
In other words, $Z=Z_3\subset \PP^4$ is a cubic threefold 
and $f(E)\in Z$ is a cA-point.
In this case, $Z^+\simeq \PP^3$, 
the map $Z\dashrightarrow Z^+$ is just a projection 
from $f(E)\in Z$, and $Y$ is a Fano threefold 
(see \cite[Table 2, $\mathrm{n^o}$ $\mathrm{15^{o}}$]{Mori1981-82}), a contradiction.
\end{proof}

\begin{case}{\bf Type \type{e_1}.}\label{case-r=2-e1}
Let $E\subset Y$ 
be the exceptional divisor and let $E^+:=\chi_*E$.
Then $Z$ is a $\QQ$-factorial Fano threefold with $\uprho(Z)=1$
and $B:=f(E)\subset Z$ is an irreducible curve.
Let $H$ be the ample generator of $\Pic(Z)$, let $H^*:=f^*H$, and $H^+:=\chi_*H^*$.
\end{case}

\begin{stheorem}{\bf Lemma.}\label{lemma-r=2-e1}
In the notation of \xref{case-r=2-e1} one of the following holds:
\begin{enumerate}
\item[\small\ref {theorem-v22-rho-2-P3}]
$Z\simeq \PP^3$, $p_a(B)=0$, $\deg B=5$;
\item[\small\ref {theorem-v22-rho-2-Q}]
$Z\simeq Q\subset \PP^4$, $p_a(B)=0$, $\deg B=5$;
\item[\small\ref {theorem-v22-rho-2-V5}]
$Z\simeq V_5\subset \PP^6$, $p_a(B)=0$, $\deg B=4$.
\end{enumerate}
\end{stheorem}
\begin{proof}
By Lemma \ref{equation-E2-5-main} we have 
\begin{equation}\label{equation-rho=2-e1}
\begin{array}{rcl}
(-K_Z)^3 &=&22+ 2K_{Y}^2 \cdot E + 2p_a(B)-2> 24,
\\[3pt]
(-K_Z)^3 &=& 22 -2K_Z \cdot B - 2p_a(B)+2.
\end{array}
\end{equation}
Hence
$\iota(Z)\ge 2$ and by Corollary \ref{Corollary-deformations-c}
we have the following possibilities:
\subsubsection*{Subcase: $\iota(Z)=4$ and $Z\simeq \PP^3$.}
Then the relations \eqref{equation-rho=2-e1} give us
\[
p_a(B)=4k,\qquad \deg B=5+k, \qquad K_{Y}^2 \cdot E =22- 4k
\]
If $k=1$, the divisor $-K_{Y}$ is in fact ample
\cite[Table 2, $\mathrm{n^o}$ $\mathrm{15^{o}}$]{Mori1981-82} and $\pi$ is an isomorphism, a contradiction.
Thus we may assume that $k\ge 2$.
Since the linear system $|-K_{Y}|$ is base point free, the curve $B$ is an intersection of quartics.
In particular, $B$ is not contained in a plane.
We get a contradiction by the Castelnuovo genus bound
(see e.g. \cite[ch. 4, Theorem 6.4]{Hartshorn-1977-ag}). 

\subsubsection*{Subcase: $\iota(Z)=3$ and $Z$ is a smooth quadric $Q\subset \PP^4$.}
As above we have 
\[
p_a(B)=3k,\qquad \deg B=5+k,\qquad K_{Y}^2 \cdot E=17-3k,\qquad 0\le k\le 5.
\]
If $k=0$, we get the case \ref{theorem-v22-rho-2-Q}.
Let $k\ge 1$. Then according to the Castelnuovo genus bound the curve 
$B$ is contained in a hyperplane.
Since $|-K_{Y}|$ is base point free, $B$ is intersection of cubics.
Hence, $\deg B\le 6$. Then 
$k=1$, $\deg B=6$, and $B$ is a complete intersection of a quadric and a cubic.
But in this case $p_a(B)=4$, a contradiction. 

\subsubsection*{Subcase: $\iota(Z)=2$ and 
$Z=Z_d\subset \PP^{d+1}$ is a del Pezzo threefold of degree $d=4$ or $5$.}
Note that $\QQ$-factorial del Pezzo threefold of degree $5$ is smooth
(and isomorphic to $V_5$, see \cite[Corollary 5.4]{Prokhorov-GFano-1}). 
As above
\[
p_a(B)=2k,\quad \deg B=2d-6+k,\quad K_{Y}^2 \cdot E =4d-10- 2k,\quad 0\le k<2d-5.
\]
If $k=0$ and $d=5$, we get the case \ref{theorem-v22-rho-2-V5}.
If $k=0$ and $d=4$, $Y$ is a Fano threefold 
as in \cite[Table 2, $\mathrm{n^o}$ $\mathrm{16^{o}}$]{Mori1981-82}, so it does not admit 
a small contraction, a contradiction.
Let $k\ge 1$. Then $p_a(B)\ge 2$ and $\deg B\ge 3$.
Note that the curve $B$ is an intersection of quadrics because 
$|-K_{Y}|$ is base point free. 
Hence $\dim \langle B\rangle\ge 4$.
Moreover, if $\dim \langle B\rangle= 4$, then 
$\deg B=5$, $d=5$, and $B=V_5\cap \langle B\rangle$
is a curve of arithmetic genus $1$ with contradicts our relation $p_a(B)=2k\ge 2$. 
Therefore, $\dim \langle B\rangle\ge 5$. 
Then, we get a contradiction by the Castelnuovo genus bound.
\end{proof}

\begin{mtheorem}{\bf Lemma.}\label{lemma-r=2-e5}
Neither $f$ nor $f^+$ is of type \type{e_5}.
\end{mtheorem}
\begin{proof}
By symmetry we may assume that $f$ is of type \type {e_5}.
Let $E\subset Y$ 
be the exceptional divisor. We have
\begin{equation}\label {equation-rho=2-e5-1}
(-K_{Y})^2\cdot E=1,
\quad (-K_{Y})\cdot E^2=-2,
\quad E^3=4.
\end{equation}
The variety $Z$ is a $\QQ$-factorial non-Gorenstein threefold with $\uprho(Z)=1$,
and $f(E)\in Z$ is a point point of type $\frac12 (1,1,1)$.
(Thus $Z$ is so-called $\QQ$-Fano threefold).
By \eqref{table-extremal-contractions} we have 
$-K_Z^3=-K_{Y}^3+1/2=45/2$.
Assume that $\Cl(Z)$ contains a torsion element, say $T$.
Since $f(E)$ is the only non-Gorenstein point of $Z$ and it is of type $\frac12 (1,1,1)$,
in its neighborhood we have $K_Z+T\sim 0$.
This means that $K_Z+T$ is a Cartier divisor. On the other hand, $(K_Z+T)^3=K_Z^3$
is not an integer, a contradiction. Therefore, the group $\Cl(Z)$ is torsion free. 
Let $A$ be the ample generator of $\Cl(Z)$.
We can write $-K_Z\sim qA$, $q\in \ZZ$. Then $45/2=-K_Z^3=q^3 A^3$.
By Proposition \ref{Proposition-factorial} the divisor $2A$ is Cartier.
Hence, $2A^3$ is an integer and so $q=1$, i.e. $K_Z$ generates 
the group $\Cl(Z)$.
Since $f_*K_{Y}=K_Z$, the group $\Pic(Y)$ is generated by $K_{Y}$ and $E$
(see \ref{length-extremal-contractions}).
Thus, for any effective divisor $D\neq E$ on $Y$, we can write 
\begin{equation*}
D\sim -\alpha K_{Y}-\beta E,\qquad \alpha\in \ZZ_{>0},\quad \beta\in \ZZ_{\ge 0}. 
\end{equation*}
Then
\begin{equation*}
(-K_{Y})\cdot D^2=22\alpha^2-2\alpha\beta-2\beta^2.
\end{equation*}
We may assume that the contraction $f^+$ is of type \type{c}, \type{d}, or
\type{e_5}. Let $D^+\subset Y^+$ be the pull-back of a line, the fiber, and the exceptional 
divisor in these cases, respectively.
Let $D:=\chi_*^{-1}D^+$. Then $(-K_{Y})\cdot D^2=2$, $0$, or $-2$, respectively.
Hence, 
\begin{equation}\label {equation-rho=2-e5-b}
11\alpha^2-\alpha\beta-\beta^2=\delta,\qquad \delta\in \{1, 0, -1\}.
\end{equation}
This equation has solutions only for $\delta=-1$, i.e. when $f^+$ is of type \type{e_5}.
But then, as above, $\Pic(Y)$ is generated by $K_{Y}$ and $D$.
Hence, $\beta=1$. In this case, \eqref{equation-rho=2-e5-b} has no solutions,
a contradiction.
\end{proof}

\begin{mtheorem}{\bf Lemma.}\label{lemma-r=2-d-d}
Both contractions $f$ and $f^+$ cannot be of type \type{d}.
\end{mtheorem}
\begin{proof}
Assume that both $f$ and $f^+$ are of type \type {d}.
Let $F$ be a fiber of $f$. Clearly,
\begin{equation*}
(-K_{Y})\cdot F^2=0,
\qquad 
(-K_{Y})^2\cdot F=K_F^2.
\end{equation*}
Let $L^+$ be a fiber of $f^+$ and let $L$ be the proper transform of $L^{+}$ under the flopping map.
Write $L\sim -\alpha K_{Y}-\beta F$, $\alpha,\beta\in \frac12 \ZZ\cup \frac13 \ZZ$
(see \ref{length-extremal-contractions}).
Then 
\begin{equation*}
\begin{array}{rclcl}
0&=& (-K_{Y})\cdot L^2 &=& 22\alpha^2-2 (K_F^2)\, \alpha\beta,
\\[3pt]
K_{F^+}^2&=& (-K_{Y})^2\cdot L &=& 22\alpha- (K_F^2)\,\beta.
\end{array}
\end{equation*}
This gives us
\begin{equation*}
11\alpha= (K_F^2)\beta= K_{F^+}^2.
\end{equation*}
But then the degree $K_{F^+}^2$ of a del Pezzo surface  $F^+$ must be divisible by 11, a contradiction.
\end{proof}

\begin{mtheorem}{\bf Lemma.}\label{lemma-r=2-c}
If both contractions $f$ and $f^+$ are not birational, then 
up to permutation we may assume that $f$ is a $\PP^1$-bundle over $\PP^2$
and $f^+$ is a del Pezzo fibration of degree $5$. 
\end{mtheorem}
\begin{proof}
By symmetry we may assume that $f$ is of type \type {c}.
Let $\Delta_f\subset Z$ be the discriminant curve of our conic bundle $f$.
Let $l\subset Z=\PP^2$ be a general line and let $F:=f^{-1}(l)$.
Then $F$ is a smooth rational surface having a conic bundle structure
with degenerate fibers over the points $l\cap \Delta_f$. Hence,
\begin{equation}
\label{equation-conic-bundle}
F^3=0, \quad K_{Y}\cdot F^2=-2,\qquad K_{Y}^2\cdot F=12 -\deg \Delta_f.
\end{equation}
We may assume that the contraction $f^+$ is of type \type{c} or \type{d}. Consider these subcases separately.

\subsubsection*{Subcase: $f^+$ is of type \type{c}}
Let $\Delta_{f^+}\subset Z^+$ be the discriminant curve of $f^+$
and let $L^+:={f^+}^{-1}(l^+)$,
where $l^+\subset Z^+=\PP^2$ is a line.
Let $L$ be the proper transform of $L^{+}$ under the flopping map.
By \ref{length-extremal-contractions} we can write 
\begin{equation*}
L\sim -\alpha K_{Y}-\beta F,\qquad \alpha,\beta\in \textstyle{\frac12} \ZZ,\qquad \alpha,\beta>0.
\end{equation*}
Similar to the proof of Lemma \ref{lemma-r=2-d-d}
we have 
\begin{equation}\label{equation-v22-rho=2-CB1}
2= (-K_{Y^+})\cdot {L^+}^2 = (-K_{Y})\cdot L^2 =
22\alpha^2-2 (12-\deg \Delta_f)\alpha\beta +2\beta^2.
\end{equation}
If $\Delta_f=\emptyset$, then the last relation can be rewritten as follows
\begin{equation*}
1+25\alpha^2= (6\alpha-\beta)^2.
\end{equation*}
It is easy to see that this equation has no solutions in $\alpha,\beta\in \frac 12 \ZZ$.
Therefore we may assume that $\Delta_f\neq\emptyset$ and $\Delta_{f^+}\neq\emptyset$.
In this case, the divisors $K_{Y}$ and $F$ generate $\Pic(Y)$.
Hence, $\alpha,\beta\in \ZZ$. By symmetry $K_{Y}$ and $L$ generate $\Pic(Y)$.
Hence, $\beta=1$. Then \eqref{equation-v22-rho=2-CB1} has the form
\begin{equation*}
12= 11\alpha+ \deg \Delta_f
\end{equation*}
Hence, $\deg \Delta_f=1$, i.e. $\Delta_f$ is a line on $Z=\PP^2$.
In this case $f^{-1}(\Delta_f)$ must be a reducible surface which contradicts
the extremal property of $f$.

\subsubsection*{Subcase: $f^+$ is of type \type{d}}
Let $L^+$ be a fiber of $f^+$ and let $L$ be the proper transform of $L^{+}$ under the flopping map.
Write $L\sim -\alpha K_{Y}-\beta F$, $\alpha,\beta\in \frac12 \ZZ$.
Then 
\begin{equation*}
\begin{array}{rcl}
0&=& (-K_{Y})\cdot L^2 = 22\alpha^2-2 (12-\deg \Delta_f)\alpha\beta +2\beta^2,
\\[3pt]
K_{L^+}^2&=& (-K_{Y})^2\cdot L = 22\alpha- (12-\deg \Delta_f)\beta.
\end{array}
\end{equation*}
The first equation 
has a rational solution only for $\deg \Delta_f=0$ and then
\begin{equation*}
K_{L^+}^2=2\beta(11\alpha/\beta- 6\beta),
\quad \alpha/\beta=1,\quad K_{L^+}^2=10\beta=5.
\end{equation*}
Since $\Delta_f=\emptyset$, the variety $Y$ is smooth and 
$f$ is a $\PP^1$-bundle (see \ref{remark-conic-bundles}). By Lemma \ref{lemma-vector-bundle}
below we get the case \ref{theorem-v22-rho-2-DP5}.
\end{proof}

\begin{mtheorem}{\bf Lemma.}\label{lemma-vector-bundle}
Let $\EEE$ be a rank-$2$ vector bundle on $\PP^2$ such that
$Y:=\PP_{\PP^2}(\EEE)$ is a generalized Fano threefold with $\g(Y)=12$. Then 
$\EEE$ is stable with $c_1(\EEE)^2+16=4c_2(\EEE)$.
\end{mtheorem}
\begin{proof}
Apply Lemma \xref{lemma-vector-bundle-2}.
It follows from  \eqref{equation-vector-bundle-K3} that $c_1(\EEE)=0$ is even and 
up to twisting $\EEE$ by a line bundle, we may assume that 
$c_1(\EEE)=0$. Moreover, $c_2(\EEE)=4$. Thus $-K_{Y}\sim 2M +3F$, where $F$ is  the pull-back 
of a line on $\PP^2$.
Since
\[
2(-K_{Y})^2\cdot M =(-K_{Y})^3- 3(-K_{Y})^2\cdot F=22-36<0,
\]
we have 
$|M|=\emptyset$ and so $H^0(\EEE)=0$.
Then $\EEE$ must be stable 
(see e.g. \cite[Ch. 4, Prop. 14]{Friedman1998}).
\end{proof} 

To finish our proof of Theorem \xref{Theorem-v22-rho-2} it remains to describe 
the right hand side of the diagram \eqref{(4.1.1)}
in cases \xref{theorem-v22-rho-2-P3}, \xref{theorem-v22-rho-2-Q}, and
\xref{theorem-v22-rho-2-V5}. 

\begin{case}\label{notation-I-II-III-a}{\bf Notation.}
Let $Z$ be a smooth Fano threefold with $\uprho(Z)=1$ and $B\subset Z$
be a smooth curve.
Let $\iota:=\iota(Z)$ be the Fano index of $Z$ and $d:=\iota-1$.
Let $H=-\frac 1\iota K_Z$ be the ample generator of $\Pic(Z)$.
Let $f: Y\to Z$ be the blowup of $B$, let $E\subset Y$ be the $f$-exceptional divisor,
and $H^*:=f^*H$.

We assume that the pair $(Z,B)$ is one of that described in Lemma \ref{lemma-r=2-e1}.
Thus $(Z,B)$ appears in the left hand side of \eqref{(4.1.1)}
in the cases \xref{theorem-v22-rho-2-P3}-\xref{theorem-v22-rho-2-V5}.
Denote by $|dH-B|$ the subsystem of $|dH|$ consisting of all 
divisors passing through $B$.
\end{case}
\begin{stheorem}{\bf Lemma.}\label{lemma-r=2-e1a}
In the above notation we have $\dim |dH-B|\ge d$.
\end{stheorem}
\begin{proof}
We have $h^0(\III_B(d))\ge h^0(\OOO_Z(d)) -h^0(\OOO_B(d))=d+1$.
\end{proof}

\begin{stheorem}{\bf Lemma.}\label{Corollary-r=2-irreducible}
If $Y$ is a generalized Fano threefold, then any member $\tilde S\in |dH^*-E|$ is irreducible.
\end{stheorem}

\begin{proof}
Let $\tilde S_1\subsetneqq \tilde S$ be an irreducible 
component other than $E$. Then $\tilde S_1\sim d' H^* -kE$, where 
$d'\le d$, $k\ge 1$ and at least one of these inequalities is strict.
Then easy computations give us $(-K_{Y})^2\cdot \tilde S_1\le 0$.
This contradicts our assumption. 
\end{proof}
\begin{stheorem}{\bf Corollary.}
In the assumptions of \xref{Corollary-r=2-irreducible}
any member $S\in |dH-B|$ is irreducible. In particular, 
$B$ is not contained in a quadric in the case $Z=\PP^3$ 
and $B$ is not contained in a hyperplane section in the case $Z=Q$. 
\end{stheorem}

\begin{case}\label{notation-I-II-III}
Now assume that $Y$ is a generalized Fano threefold, i.e. the divisor
$-K_Y$ is nef and 
the morphism defined by the linear system  $|-nK_Y|$
does not contract divisors.
Let $\chi: Y \dashrightarrow Y^+$ be the flop, let
$H^+:=\chi_*H^*$, and $E^+:=\chi_*E$. The description of the right hand side of the diagram \eqref{(4.1.1)}
in cases \xref{theorem-v22-rho-2-P3}, \xref{theorem-v22-rho-2-Q}, and
\xref{theorem-v22-rho-2-V5} is an immediate consequence of the following.
\end{case}

\begin{mtheorem}{\bf Proposition.}\label{lemma-r=2-e1b}
In the above notation, $dH^+-E^+$ is a supporting divisor for the contraction 
$f^+$, i.e. the morphism $f^+$is given by the linear system $|n(dH^+-E^+)|$, $n\gg 0$.
\end{mtheorem}
\begin{proof}
Let $\tilde S\in |dH^*-E|$ and $\tilde S^+:=\chi_*\tilde S$.
Since $-\tilde S\sim K_{Y}+H^*$ is $\pi$-ample, $\tilde S^+$ must be $\pi^+$-ample.
Since the linear system $|\tilde S^+|$ has no fixed components and the contraction $f^+$ is not small,
$\tilde S^+$ is non-negative on the fibers of $f^+$.
Therefore, $\tilde S^+$ is nef.
It remains to show that $\tilde S^+$ is not ample. Assume the converse.
Let $l$ be a curve in a fiber of $f^+$ that attains the minimum of 
$\mu_{f^+}$ (see \eqref{equation-length-extremal-contractions}).
Since $|H^+|$ has no fixed components and $H^+$ is negative on flipped curves, we have $H^+\cdot L>0$.
Then $\mu_{f^+}=-K_{Y^+}\cdot l=\tilde S^+\cdot l+H^+\cdot l\ge 2$.
Moreover, if $\mu_{f^+}=2$, then $\tilde S^+\cdot l=H^+\cdot l=1$.
In this case the supporting divisor for $f^+$ has the form $\tilde S^+-H^+$ 
and so $|\tilde S-H^*|\neq \emptyset$. 
On the other hand, 
$(-K_{Y})^2\cdot (\tilde S-H^*)\le 0$.
This contradicts our assumption that $Y$ is a generalized Fano threefold.
Hence, $\mu_{f^+}=3$. 
Then $Z^+\simeq \PP^1$ and $f^+$ is a $\PP^2$-bundle.
As above, the only possibility is the following:
$\tilde S^+\cdot l=1$, $H^+\cdot l=1$, and the supporting divisor of $\varphi^+$ is 
$2\tilde S^+-H^+$. Since $\dim Z^+=1$, $(2\tilde S^+-H^+)^2=0$.
On the other hand, $(2\tilde S^+-H^+)^2\cdot K_{Y^+}=(2\tilde S-H^*)^2\cdot K_{Y}<0$
by \eqref{equation-E1-main}, a contradiction.
\end{proof}

This concludes our proof of 
Theorem \xref{Theorem-v22-rho-2}.

\section{Fano threefolds of types \xref{theorem-v22-rho-2-P3}, \xref{theorem-v22-rho-2-Q}, and
\xref{theorem-v22-rho-2-V5}.}
\label{section-I-II-III}
First we show that the possibilities 
\xref{theorem-v22-rho-2-P3}, \xref{theorem-v22-rho-2-Q} and
\xref{theorem-v22-rho-2-V5} of Theorem \xref{Theorem-v22-rho-2}
really occur.
More precisely, we show that any pair $(Z,B)$
satisfying the conditions of 
\xref{theorem-v22-rho-2-P3}-\xref{theorem-v22-rho-2-V5} of Theorem \xref{Theorem-v22-rho-2}
generates diagram \eqref{(4.1.1)}.
\begin{case}\label{notation-I-II-III-b}
Let $(Z,B)$ be as in \ref{notation-I-II-III-a}.
By Lemma \ref {lemma-r=2-e1a}
$\dim |dH-B|>0$. Assume additionally that $B$ is not contained in a quadric in the case $Z=\PP^3$ 
and $B$ is not contained in a hyperplane section in the case $Z=Q$. 
This implies that any member $S\in |dH-B|$ is irreducible.
\end{case}

\begin{stheorem}{\bf Lemma.}\label{Corollary-r=2-e1a}
Any member $S\in |dH-B|$ is del Pezzo surface of degree $d$ with at worst Du Val singularities.
\end{stheorem}

\begin{proof}
First we claim that $S$ is normal.
Assume that $S$ is singular along a curve $J$. Then 
$J$ is contained in $\Bs |dH-B|$ by Bertini's theorem. 
The intersection of two members $S, S'\in |(\iota-1) H -B|$ contains
$J$ with multiplicity $\ge 4$
and the curve $B$. Then $(\iota-1)^2H^3=S\cdot S'\cdot H\ge H\cdot B +4H\cdot B $.
In the cases \xref{theorem-v22-rho-2-Q} and
\xref{theorem-v22-rho-2-V5} this gives a contradiction.
In the case \xref{theorem-v22-rho-2-P3} we have an equality
$S\cdot S'= B +4 J$. Since $\dim |(\iota-1) H -B|=3$ in this case,
we can take $S$ and $S'$ so that they have a common point $P\notin B\cap J$.
Then $S\cdot S'\supsetneqq B +4 J$. This means that $S$ and $S'$ have a common component
and $S$ is reducible. The contradiction proves our claim.

Further, by the adjunction formula $-K_{S}=H|_{S}$ is ample, i.e. $S$ is a del Pezzo surface.
If the singularities of $S$ are worse than Du Val, then it is a cone over an elliptic curve 
\cite{Hidaka1981}. But this is impossible because $S$ contains a rational curve $B$
of degree $>1$.
\end{proof}
\begin{stheorem}{\bf Corollary.}
In the above notation any member $\tilde S\in |dH^*-E|$ 
is a del Pezzo surface of degree $d$ with at worst Du Val singularities
and the restriction map $f_S: \tilde S\to f(\tilde S)=S$ is crepant.
\end{stheorem}
\begin{proof}
 By the adjunction formula $K_{\tilde S}=-f^*H|_{\tilde S}=f^* K_S$. 
\end{proof}

\begin{mtheorem}{\bf Proposition.}\label{Proposition-LemmaI-II-IIIa}
Notation as in \xref{notation-I-II-III-a}.
Let $S\in |dH-B|$ be a general member and let $\tilde S\subset Y$ 
be the proper transform of $S$.
Then $S$ is smooth and $f_S: \tilde S\to S$ is an isomorphism.
The blowup $f: Y\to Z$ can be completed to the diagram \eqref{(4.1.1)}.
The $\pi$-exceptional locus consists of exactly one 
smooth rational curve $\Upsilon$ which is a $(-1)$-curve on 
$\tilde S$ and
$\NNN_{\Upsilon/ Y}\simeq\OOO_{\PP^1}(-1)\oplus \OOO_{\PP^1}(-1)$.
\end{mtheorem}

\begin{proof}
We claim that $Y$ is a generalized Fano threefold.
Assume that $K_Y\cdot \Upsilon\ge 0$ for some irreducible curve $\Upsilon$.
Clearly, $\Upsilon$ is not a fiber of the ruling $E\to B$.
Hence, $H^*\cdot \Upsilon>0$ and $\tilde S\cdot \Upsilon<0$. 
In particular, $\Upsilon\subset \tilde S$, $\Upsilon\subset \Bs|dH^*-E|$,
and there are at most a finite number of such curves. 
On the minimal resolution $S_{\min}$ of $\tilde S$, the Mori cone 
is generated by $(-1)$- and $(-2)$-curves 
(because $-K_{S_{\min}}$ is nef and big).
Hence we may assume that $-K_{\tilde S}\cdot \Upsilon=1$
and so $H\cdot f(\Upsilon)=-K_S\cdot f(\Upsilon)=1$, i.e. $f(\Upsilon)$ is a line on $Z$.
Since $K_Y\cdot \Upsilon\ge 0$, we have $E\cdot \Upsilon\ge d+1$
and so $f(\Upsilon)$ is a $(d+1)$-secant line to $B\subset S\subset \PP^d$. 
It is easy to see that $B$ has no $(d+2)$-secant lines.
Therefore, $E\cdot \Upsilon= d+1$ and $K_Y\cdot \Upsilon= 0$, i.e. $-K_Y$ is nef.
Obviously, $-K_Y$ is big.
This proves our claim.

Thus the blowup $f: Y\to Z$ can be completed to the diagram \eqref{(4.1.1)}
and by Theorem \xref{Theorem-v22-rho-2} we have cases 
\xref{theorem-v22-rho-2-P3}, \xref{theorem-v22-rho-2-Q}, 
\xref{theorem-v22-rho-2-V5}.
Let $\tilde S^+=\chi_*\tilde S\subset Y^+$ 
be the proper transform of $\tilde S$.
Then 
$\tilde S^+$ is a smooth del Pezzo surface of degree $d+1$.
Since the linear system $|\tilde S^+|$ is base point free, $\tilde S^+$ meets flopped curves 
$\Upsilon _i^+$
transversely. Hence $\tilde S$ contains all the flopping curves $\Upsilon _i$.
By the Zariski main theorem $\chi$ induces a morphism $\chi_S: \tilde S\to \tilde S^+$
that contracts $\cup \Upsilon _i$. 
By the Noether formula 
$\uprho(\tilde S)\le 10-d$ and $\uprho(\tilde S^+)=9-d$.
Since the exists at least one flopping curve, 
$\uprho(\tilde S)=\uprho(\tilde S^+)+1$ and $\uprho(\tilde S)= 10-d$.
Hence $\tilde S$ is smooth and $\chi_S: \tilde S\to \tilde S^+$
is a blowup of a single point.
Thus the $\pi$-exceptional locus consists of exactly one 
smooth rational curve $\Upsilon=\Upsilon_1$ which is a $(-1)$-curve on 
$\tilde S$. Then 
$\NNN_{\Upsilon/ Y}$ contains a subbundle
$\NNN_{\Upsilon/ \tilde S}\simeq \OOO_{\PP_1}(-1)$.
Since $\deg \NNN_{\Upsilon/ Y}=-2$, we have the desired splitting.
If the restriction $f_S: \tilde S\to S$ is not an isomorphism,
then the intersection $\tilde S\cap E$ is reducible: it contains 
a section and a fiber of the ruling $E\to B$.
But this contradicts Bertini's theorem applied to $|\tilde S|\bigr|_{E}$.
\end{proof}

\begin{stheorem}{\bf Corollary.}
In the cases \xref{theorem-v22-rho-2-P3}, \xref{theorem-v22-rho-2-Q},
\xref{theorem-v22-rho-2-V5} of Theorem \xref{Theorem-v22-rho-2} the 
singular locus of $X$ consists of exactly one ordinary double point.
\end{stheorem}

\begin{case}{\bf Construction.}\label{Proposition-I-II-III-inverse}
Let $Z$ be either $\PP^3$, a (possibly singular) quadric in $\PP^4$ or a (possibly singular)
quintic del Pezzo 
threefold in $\PP^6$. Put $d:=\iota(Z)-1$.
Let $H$ be the ample generator of $\Pic(Z)$ and let $S=S_d \in |dH|$
be a smooth member. Thus 
$S=S_d\subset \PP^d$ is a smooth del Pezzo surface of degree $d=3$, $4$, or $5$.
Regard $S$ as the blowup $\sigma: S\to \PP^2$ of $9-d$ points in general position.
Let $\mathbf h$ be the class of line on $\PP^2$,
$\mathbf h^*:=\sigma^*\mathbf h$ and let $\mathbf e_1$, \dots, 
$\mathbf e_{9-d}$ be the $\sigma$-exceptional divisors.
Let $B\subset S$ be a reduced curve such that 
$B\sim 2\mathbf h^*-\mathbf e_1$
in the cases $d=3$, $4$ and $B\sim 2\mathbf h^*-\mathbf e_1-\mathbf e_2$
in the case $d=5$. 
Then $p_a(B)=0$. Moreover, $-K_S\cdot B=5$ in the cases $d=3$, $4$ and
$-K_S\cdot B=4$ 
in the case $d=5$.

Let $f: Y\to Z$ be the blowup of $B$ and let $\tilde S\subset Y$ be the proper transform of $S$.
Then the singularities of $Y$ are 
at worst terminal Gorenstein (cf. Theorem \ref{theorem-Cutkosky}\type{e_1}). Clearly,
$\tilde S\simeq S$.
Up to this identification we can write
$-K_{Y}|_{\tilde S}\sim -(d+1)K_S-B$.
It is easy to check that this divisor is nef and big.
Moreover, $-K_{Y}|_{\tilde S}$ has positive intersection number 
with all curves except for one, which is a line on $S$. 
Denote it by $\Lambda$.
If $-K_{Y}\cdot \Upsilon \le 0$ for some curve $\Upsilon$, then
$\tilde S\cdot \Upsilon <0$. Hence $\Upsilon \subset \tilde S$ 
and $\Upsilon =\Lambda$. This shows that 
the divisor $-K_{Y}$ is nef and $\Lambda$ is the only $K$-trivial curve on $Y$. 
By using \eqref{equation-E1-main} one can compute 
$-K_Y^3=22$.
So, $Y$ is a generalized Fano threefold
of genus $12$.
\end{case}
Thus the possibilities \xref{theorem-v22-rho-2-P3}, \xref{theorem-v22-rho-2-Q} ╔
\xref{theorem-v22-rho-2-V5} of Theorem \xref{Theorem-v22-rho-2} 
occurs.

In the above construction the curve $B$ can be reducible and then $\rr(Y)-\rr(Z)$
equals to the number of component of $B$. 
\begin{scase}{\bf Example.}\label{example-rho=8}
Let $Z=V_5'$ be the the quintic del Pezzo threefold with $\rr(Z)=4$
(see \cite{Prokhorov-GFano-1}) and let $B$ be a combinatorial
chain of four lines contained in a smooth hyperplane section $S$. 
Then the above construction give us a generalized Fano threefold $Y$ with 
$\rr(Y)=8$. 
\end{scase}
Another application of our construction \ref{Proposition-I-II-III-inverse}
is the following.

\begin{mtheorem}{\bf Theorem (cf. \cite{Furushim2006}).}\label{compactification}
There exists a Fano threefold $X=X_{22}\subset \PP^{13}$
of type \xref{theorem-v22-rho-2-V5} and a reducible hyperplane
section $A=A_1\cup A_2$ such that $X\setminus A\simeq \mathbb{A}^3$. 
\end{mtheorem}

\begin{proof}
Recall that the normal bundle of a line on $V_5$ has 
the form $\OOO_{\PP^1}(-a)\oplus \OOO_{\PP^1}(a)$ with $a=0$ or $1$
(see e.g. \cite[Lemma 4.2.1]{Iskovskikh-Prokhorov-1999} or \cite{Furushima1989a}).
Let $\Lambda\subset V_5$ be a line with 
$\NNN_{\Lambda/V_5}\simeq \OOO_{\PP^1}(-1)\oplus \OOO_{\PP^1}(1)$
and let $R$ be the ruled surface swept out by lines meeting $\Lambda$.
Then $R$ is a hyperplane section of $V_5$ and its normalization $R'$ is isomorphic to 
$\FF_3$ \cite{Furushima1989a}. 
Moreover, $V_5\setminus R\simeq \mathbb{A}^3$ (loc. sit.).
The map $\nu: R'\to R$ is an isomorphism outside of $\Lambda$ and $\nu^{-1}(\Lambda)$ is the 
union of the negative section $\Sigma\subset R'=\FF_3$ and 
a fiber $\Gamma_1$. 
Further, the map
$\nu : R'\to R\subset \PP^5$ is given by a codimension $1$ subsystem of the linear system $|\Sigma+4\Gamma_1|$.
Let $B'\subset R'$ be a general member of $|\Sigma+3\Gamma_1|$.
Then $B:=\nu(B')$ is a smooth rational curve of degree $4$
and the line $\Upsilon:=\nu(\Gamma_0)$ 
is a $2$-secant of $B$ for some fiber $\Gamma_0\subset R'$.
Thus $V_5\cap \langle B\rangle=B\cup \Upsilon$.
Then we can apply Proposition \ref{Proposition-LemmaI-II-IIIa}
and get a
Fano threefold $X_{22}$. By our construction $X_{22}\setminus (\pi(f^{-1}(R)))\simeq V_5\setminus R\simeq \mathbb{A}^3$.
\end{proof}

\section{Fano threefolds of type \xref{theorem-v22-rho-2-DP5}}
\label{section-IV}
In this section we investigate the case \xref{theorem-v22-rho-2-DP5} of Theorem \xref{Theorem-v22-rho-2}.
The relation between singular Fano threefolds $X_{22}$ and 
certain rank-$2$ vector bundles 
was noticed by Mukai \cite[Remark 5]{Mukai-1989}.
This observation is based on the explicit description of stable rank-$2$ vector bundles
on $\PP^2$ with even $c_1$ \cite{Barth1977}.

\begin{case}
Let $\EEE$ be a stable rank-$2$ vector bundle on $\PP^2$ with $c_1(\EEE)=0$, $c_2(\EEE)=4$
\cite{Barth1977}.
Let $Y:=\PP_{\PP^2}(\EEE)$ and let $f:Y \to \PP^2$ be the projection.
Let $F$ be the pull-back of a line $l\subset \PP^2$ and let $M$ be the tautological divisor.
Then
\begin{equation}
\label{equation-vector-bundle-Hirsch-a}
M^3=-4, \qquad 
M^2\cdot F=0, \qquad M\cdot F^2=1, \qquad 
F^3=0
\end{equation}
(see \eqref{equation-vector-bundle-Hirsch}).
For $n\le 2$, by the Serre duality we have 
\begin{equation*}
H^2(\EEE(n))=
H^0\left(\EEE(n)\otimes \det \EEE(n)^{\vee}\otimes \upomega_{\PP^2}\right)^{\vee}
=
H^0(\EEE(n-3))^{\vee}=0.
\end{equation*}
Then, by the Riemann-Roch theorem 
\begin{equation*}
h^0\left(\EEE(n)\right) 
\ge \frac 12\left(c_1(\EEE(n))^2 - 2 c_2(\EEE(n)) - K_{\PP^2} \cdot c_1(\EEE(n))\right) +2=n^2 + 3n - 2.
\end{equation*}
In particular,
\begin{equation}
\begin{array}{rcl}
\dim |M+F|&=&h^0(\EEE(1))-1\ge 1,
\\[4pt]
\dim |M+2F|&=&h^0(\EEE(2))-1\ge 7.
\end{array}
\end{equation}
Since $\EEE$ is stable, $H^0(\EEE)=0$ (see e.g. \cite[Ch. 4, Prop. 14]{Friedman1998})
and so $|M|=\emptyset$.
Hence any member $S\in |M+F|$ is irreducible.
Let $S,\, S'\in |M+F|$ be general members and let $\Gamma:= S\cap S'$
(scheme-theoretically). 
\begin{stheorem}{\bf Claim.}\label{lemma-vector-bundle-1}
In the above notation, we have 
$-K_Y\cdot \Gamma=0$, in particular, $Y$ is not a Fano threefold.
Moreover, 
the image $f(\Gamma)$ is a \textup(possibly degenerate\textup) conic.
\end{stheorem}
\begin{proof}
Using \eqref{equation-vector-bundle-Hirsch-a} we obtain
$(M+F)^2\cdot K_Y=0$ and $(M+F)^2\cdot F=2$.
\end{proof}

\begin{mtheorem}{\bf Proposition.}\label{Proposition-vector-bundle}
In the above notation the following conditions are equivalent:
\begin{enumerate}
\item \label{Proposition-vector-bundle-1}
$Y$ is a generalized Fano threefold,
\item \label{Proposition-vector-bundle-2}
$\EEE(2)$ is ample \textup(or, equivalently, $M+2F$ is ample\textup),
\item\label{Proposition-vector-bundle-3}
$\Gamma$ is reduced and irreducible,
\item\label{Proposition-vector-bundle-3a}
$\Gamma$ is a smooth rational curve,
\item\label{Proposition-vector-bundle-4}
a general member $S\in |M+F|$ is a smooth quartic del Pezzo surface.
\end{enumerate}
\end{mtheorem}

\begin{proof}
\ref{Proposition-vector-bundle-1} $\Longrightarrow$ \ref{Proposition-vector-bundle-2}
Since $\uprho(Y)=2$, the Mori cone $\overline{\operatorname{NE}}(Y)$ is generated by two extremal rays
$\mathcal R_0$ and $\mathcal R$.
We may assume that $\mathcal R_0$ is generated by curves in the fibers of the projection $f: Y\to \PP^2$.
Clearly, all effective divisors are non-negative on $\mathcal R_0$.
Since $Y$ is a generalized Fano and $-K_Y\cdot \Gamma=0$, the ray $\mathcal R$ is generated by 
the class of $\Gamma$. 
Hence $(M+2F)\cdot \mathcal R> 0$.
By the Kleiman ampleness criterion, $M+2F$ is ample and so $\EEE(2)$ is.

\ref{Proposition-vector-bundle-2} $\Longrightarrow$ \ref{Proposition-vector-bundle-3}
By the assumption $M+2F$ is ample. Since $(M+2F)\cdot \Gamma=(M+2F)\cdot (M+F)^2=1$, $\Gamma$ is reduced and irreducible.

\ref{Proposition-vector-bundle-3} $\Longrightarrow$ \ref{Proposition-vector-bundle-1}
Suppose that $-K_Y\cdot C\le 0$ for some irreducible curve $C$.
Then $(M+F)\cdot C<0$. Hence $C$ is contained in the base locus of $|M+F|$.
Since $\Gamma$ is reduced and irreducible, $C=\Gamma$.
In this case, $-K_Y\cdot C=0$ and so the divisor $-K_Y$ is nef.
Since $-K_Y^3=22$, it is big. 
Moreover, $\Gamma$ is the only curve contracted by $|-K_Y|$.

\ref{Proposition-vector-bundle-3} $\Longrightarrow$ \ref{Proposition-vector-bundle-3a}
follows by the adjunction formula, \ref{Proposition-vector-bundle-3a} $\Longrightarrow$ \ref{Proposition-vector-bundle-3}
is obvious, and the implication \ref{Proposition-vector-bundle-3a} \& \ref{Proposition-vector-bundle-2} $\Longrightarrow$ \ref{Proposition-vector-bundle-4}
follows by Bertini's theorem and adjunction.
Let us show \ref{Proposition-vector-bundle-4} $\Rightarrow$ \ref{Proposition-vector-bundle-3}.
By adjunction on $S$ we have 
\begin{equation*}
-K_S\cdot \Gamma=(M+2H)\cdot (M+H)^2 =1.
\end{equation*}
Since $-K_S$ is ample, $\Gamma$ is reduced and irreducible.

The assertion is proved.
\end{proof}
\end{case}

\begin{case}
From now on we assume that $\EEE$ satisfies the equivalent conditions of Proposition 
\ref{Proposition-vector-bundle}.
Let $\pi: Y\to X=X_{22}\subset \PP^{13}$ be the anticanonical map.
It follows from the proof of \ref{Proposition-vector-bundle} that 
$\Gamma$ is the only curve contracted by $|-K_Y|$ (i.e. by the corresponding morphism).
Therefore, the singular locus of $X$ consists of one point 
$P:=\pi(\Gamma)$. Let $\Gamma^+\subset Y^+$ be the flopped curve.
Thus $\pi^+(\Gamma^+)=P$.
Then one can restore the diagram \eqref{(4.1.1)} with the map 
$f^+$ being a degree $5$ del Pezzo fibration by Theorem \ref{Theorem-v22-rho-2}.
The right hand side of \eqref{(4.1.1)} is explicitly described by
K. Takeuchi \cite{Takeuchi-2009}.
\end{case}

\begin{scase}{\bf Remarks.}
\begin{enumerate}
\item By Claim \ref{lemma-vector-bundle-1} and Proposition \ref{Proposition-vector-bundle}, \ref{Proposition-vector-bundle-3a}
the image $\Omega:=f(\Gamma)$ is a smooth conic. 
\item 
Since $\NNN_{\Gamma/Y}$ contains a subbundle $\NNN_{\Gamma/S}\simeq \OOO_{\PP^1}(-1)$,
we have $\NNN_{\Gamma/Y}\simeq \OOO_{\PP^1}(-1)\oplus \OOO_{\PP^1}(-1)$.
In particular, the image $\pi(\Gamma)\in X$ is an ordinary double point.
\item 
By the Kawamata-Viehweg vanishing 
theorem $H^i(\OOO_Y(M+F))=0$ for $i>0$. Hence, $\dim |M+F|=1$.
\item 
Since the base locus of the pencil $|M+F|$ is a smooth curve $\Gamma$,
\emph{any} member $S\in |M+F|$ has at worst isolated singularities
and $\Sing(S)\cap \Gamma=\emptyset$. Moreover, $\Gamma$ is a $(-1)$-curve on $S$
and the singularities of $S$ are at worst Du Val \cite{Hidaka1981}.
\end{enumerate}
\end{scase}

\begin{case}
The map $f_S: S\to \PP^2$ is the blowup of the ideal sheaf of some 
zero-dimensional scheme $\Xi\subset \Omega\subset \PP^2$ of length 
$5$. 
If $S\in |M+F|$ is general, then $S$ is a smooth del Pezzo surface of degree $4$
and 
$f_S: S\to \PP^2$ is the blowup of five points in general position.
\end{case}

\begin{scase}
Denote $T:=f^{-1}(\Omega)$.
It is easy to show that $T\simeq \FF_6$ and $\Gamma\subset T$ is the negative section.
Since $(M+F)^2\cdot T=4$, any member of the restriction $|M+F|\bigr|_{T}$ has the form 
$\Gamma + \Upsilon_1+\cdots+\Upsilon_5$, where $\Upsilon_i$ are fibers.
In other words, $|M+F|\bigr|_{T}=\Gamma + f_T^* \mathfrak g$, where $\mathfrak g=\mathfrak g^1_5$ is 
a pencil (linear series) on $\Omega$ of degree $5$. Thus $\EEE$ defines a $\mathfrak g^1_5$ on 
$\Omega\simeq \PP^1$. 
\end{scase}

\begin{case}
In general, a stable rank-$2$ vector bundle $\EEE$ on $\PP^2$ with $c_1(\EEE)=0$, $c_2(\EEE)=4$
suits to the following exact (non-unique) sequence
\begin{equation}\label{equation-exact-sequence}
0 \longrightarrow \OOO_{\PP^2} \longrightarrow \EEE(1) \longrightarrow 
\OOO_{\PP^2}(2)\otimes \III_{\Xi} \longrightarrow 0,
\end{equation}
where $\III_{\Xi}$ 
is the ideal sheaf of a zero-dimensional subscheme $\Xi\subset \PP^2$ of length 
$5$ which is not contained in a line (see e.g. \cite[Ch. 4, Prop. 14]{Friedman1998}). 
Conversely, for $\Xi=\{P_1,\dots,P_5\}$ the extension \eqref{equation-exact-sequence}
corresponds to an element $s$ of the $5$-dimensional vector space
\[
\operatorname{Ext}^1_{}(\III_{\Xi}(2), \OOO_{\PP^2})\simeq \OOO_\Xi(1)=\bigoplus \OOO_{P_i}(1).
\]
and for general choice of $s$ the sheaf $\EEE(1)$ is locally free \cite[\S 5]{Barth1977}.
This shows that for a general choice of $\EEE$ (in the sense of moduli)
the vector bundle $\EEE$ satisfies the conditions of \ref{Proposition-vector-bundle}.
\end{case}

\begin{case}
We say that an irreducible curve $\Lambda$ on a generalized Fano 
threefold $U$ is a \emph{line} if $-K_U\cdot \Lambda=1$. 
Denote by $\mathfrak L(U)$ the union of those components the Hilbert scheme
$\operatorname{Hilb}(U)$
that contain classes of lines.
Roughly speaking, we may regard $\mathfrak L(U)$ as the scheme parametrizing lines on $U$.
\end{case}

\begin{stheorem}{\bf Lemma.}\label{lemma-lines-v22-IV}
Let $\Lambda\subset Y$ be a line. 
Then $F\cdot \Lambda=1$ and $\Gamma \cap \Lambda=\emptyset$. 
\end{stheorem}
\begin{proof}
Since $M+2F$ is ample, we have $(M+2F)\cdot \Lambda=1$ and $(M+F)\cdot \Lambda=0$.
Hence $\Lambda$ is disjoint from a general member of the pencil $|M+F|$.
\end{proof}
\begin{stheorem}{\bf Corollary.}\label{Corollary-lines-v22-IV}
Any line on $X$ is contained in a smooth locus and so $\mathfrak L(Y)\simeq\mathfrak L(X)\simeq\mathfrak L(Y^+)$. 
In particular, $\mathfrak L(X)$ parametrizes lines in the fibers of $f^+$.
\end{stheorem}
\begin{stheorem}{\bf Corollary.}
For a general line $\Lambda\subset X$ there are exactly three other lines meeting $\Lambda$. 
\end{stheorem}
Note that the same holds for \textit{general} non-singular Fano thereefold of genus 12
(see \cite[Proof of Theorem 6.1]{Iskovskih1978a}).

\begin{case}
Now we give another description of varieties $X_{22}$ of type \xref{theorem-v22-rho-2-DP5} which does not use
the vector bundle techniques.
\end{case}

\begin{stheorem}{\bf Lemma.}
Let $g: \hat X\to X$ be the blowup of the singular point $P\in X$.
Then $\hat X$ is a \textup(smooth\textup) Fano threefold with $\uprho(\hat X)=3$
and $(-K_{\hat X})^3=20$.
\end{stheorem}
\begin{proof}
By Corollary \ref{Corollary-lines-v22-IV} the variety $X=X_{22}\subset \PP^{13}$ 
contains no lines passing through $P$ and so $T_{P,X}\cap X=\{P\}$.
\end{proof}

\begin{scase}
Clearly, the exceptional divisor $R\subset \hat X$ is isomorphic to $\PP^1\times \PP^1$ 
and $\hat X$ admits two contractions $\tau:\hat X\to Y$ and $\tau^+: \hat X\to Y^+$ 
such that $\tau(R)=\Gamma$ and $\tau^+(R)=\Gamma^+$.
Thus we have the product map 
\[
h=(f\comp \tau) \times (f^+\comp \tau^+) :\hat X\to \PP^2\times \PP^1.
\]
We get the following extension of the diagram \eqref{(4.1.1)}:
\begin{equation}\label{equation-last-diagram}
\vcenter{
\xymatrix{
&&\PP^2\times \PP^1\ar@/^2pc/[ddrr]^{p_1}\ar@/_2pc/[ddll]_{p_2}&&
\\
&Y\ar[dl]^f\ar[dr]_{\pi}&\hat X\ar[r]^{\tau^+}\ar[l]^{\tau}\ar[d]^{g}\ar[u]^{h}&Y^+\ar[dl]^{\pi^+}\ar[dr]_{f^+}
\\
\PP^2&&X&&\PP^1
}
} 
\end{equation}
Since $\hat X$ is a Fano threefold, $h$ is an extremal contraction.
Since $\bb_3(\hat X)=0$ and $(-K_{\PP^2\times \PP^1})^3-(-K_{\hat X})^3=34$,
the contraction $h$ is the blowup of a smooth rational curve $C\subset \PP^2\times \PP^1$
with $-K_{\PP^2\times \PP^1}\cdot C=16$ (see \eqref{equation-E1-main}).
Then it is easy to see that $C$ is a curve of bidegree $(2,5)$ 
such that the restriction $p_2 |_C$ is an embedding (see \cite[Table 3, $\mathrm{n^o}$ $\mathrm{5^{o}}$]{Mori1981-82}). Clearly,
$p_2(C)=\Omega$ and $h(R)=p_2^{-1}(\Omega)$.
The $h$-exceptional divisor is the proper transform of $f^{-1}(\Omega)$.
Let $G:=p_2^{-1}(\Omega)=h(R)$. Then $G\simeq \PP^1\times \PP^1$ and 
$C\subset G$ is a curve of bidegree $(5,1)$. The projection $C\to\PP^1$ to the first factor 
is defined by a linear system $\mathfrak g^1_5$.
Conversely, a $\mathfrak g^1_5$ on $\PP^1$ defines an embedding $\PP^1\hookrightarrow \PP^1\times \PP^1=G$
as a divisor of bidegree $(5,1)$ and therefore it defines the whole diagram 
\eqref{equation-last-diagram}, where $X$ is a Fano threefold with $\g(X)=12$. 
\end{scase}

\input appendix.tex

\section{Proof of Theorem \ref{Theorem-main}}
To show that in the notation of Theorem \ref{Theorem-main} the case $\rr(X)>2$ 
does not occur, we need some upper bounds of the rank $\rr(X)$.
We start with the following result.

\begin{mtheorem}{\bf Proposition.}\label{Proposition-le10}
Let $X=X_{22}\subset\PP^{13}$ be a Fano threefold of the main series with $\g(X)=12$.
Assume that $X$ does not contain planes.
Then $\rr(X)\le 10$.
\end{mtheorem}
Note that proofs of similar estimates of $\rr(X)$ in \cite[Sect. 3.2]{Kaloghiros2011} 
contain some gaps (see footnote (${}^{\ref{footnote-Kaloghiros2011}}$) on page \pageref{footnote-Kaloghiros2011}).
\begin{proof}
Assume that $\rr(X)\ge 11$. Let $\pi:\tilde X\to X$ be a factorialization.
Run the $K$-MMP on $\tilde X$. We get diagram \eqref{equation-MMP}, where
all the maps $\varphi_i$ are divisorial contractions.

\begin{stheorem}{\bf Claim.}\label{Lemma-degree-3}\label{Lemma-degree-2}
$\varphi$ contracts a surface of degree $2$.
\end{stheorem}
\begin{proof}
Let the minimal degree of a surface contracted by $\varphi$ equals $d$.
Assume that $d\ge 3$. 
We also may assume that $\tilde X_N$ does not admit birational contractions 
and $\uprho(\tilde X_N)\le 3$
by Lemma \ref{Lemma-surface-rho=2}.
Hence, $N\ge \rr(X)-3\ge 8$.
By \eqref{equation-E1-main-0}
\begin{equation}\label{equation-degree-last}
-K_{\tilde X_N}^3\ge 22+(2d-2)N\ge 54.
\end{equation}
If $-K_{\tilde X_N}^3=54$, then by Lemma \ref{Proposition-XN} 
we have $\uprho(\tilde X_N)\le 2$ and so $N\ge 9$. In this case, $-K_{\tilde X_N}^3\ge 22+4N\ge 56$,
a contradiction. Therefore, $-K_{\tilde X_N}^3>54$ and again by Lemma \ref{Proposition-XN} 
we have only one possibility:
$\tilde X_N\simeq \PP^3$. Then $N=10$ and $d=3$. 
Moreover, similar to \eqref{equation-degree-last} we have $-K_{\tilde X_{N-1}}^3\ge 58$.
By the classification \cite{Mori1981-82} the only possibility for the smoothing of $X_{N-1}$ is 
\cite[Table 2, $\mathrm{n^o}$ $\mathrm{36^o}$]{Mori1981-82}.
But in this case $X_{N-1}$ must contain a plane, a contradiction.
\end{proof}
Let $S\subset \tilde X$ be an irreducible surface of degree $2$ contracted by $\varphi$.
Run the $S$-MMP on $\tilde X$. 
We claim that $S$ is contracted during this process.
Indeed, otherwise 
after a number of flops we get a model $\tilde X'$ such that 
the proper transform $S'\subset \tilde X'$ of $S$ is nef. 
Since $\tilde X'$ is a generalized Fano threefold,
some multiple of $S'$ must be movable. Hence the image $\bar S\subset X$ of $S'$ under the anticanonical 
morphism is movable either. On the other hand, $\bar S$ is a quadric and the intersection 
$\bar S\cap H$ with a general member $H\in |-K_X|$ is a smooth rational curve.
Since $H$ is a smooth K3 surface, the curve $\bar S\cap H$ cannot be movable.
The contradiction shows that after a number of flops we must contract the proper transform of
$S$. Replacing $\tilde X$ with another factorialization, we may assume that
$\varphi_1$ contracts a surface $E_0=S$ of degree $2$.
We claim that $B_1=\varphi_1(E_0)$ is either a point or a smooth rational curve.
Indeed, otherwise $B_1$ is a rational curve which is singular at some point $P$.
Since $\varphi_1$ is the blowup of a locally planar curve $B_1$ contained in the smooth locus of $\tilde X_1$
(see Theorem \ref{theorem-Cutkosky}), 
easy local computations show that the exceptional divisor $E_1$ is singular along 
the fiber $C:=\varphi_1^{-1}(P)$.
Since the quadric $\pi(E_1)$ is normal, $C$ must be contracted by $\pi$.
This contradicts the fact that $\varphi_1$ is a Mori contraction and proves our claim.
Thus by \eqref{table-extremal-contractions} and \eqref {equation-E1-main-0}
$-K_{\tilde X}^3=24$.
Taking Proposition \ref{proposition-24} into account we obtain 
\begin{equation*}
\rr(X)=\uprho(\tilde X)=\uprho(\tilde X_1)+1=\rr(X_1)+1\le 10. 
\end{equation*}
This proves the proposition.
\end{proof}

\begin{mtheorem}{\bf Lemma.}\label{Lemma-surface-degree=10}
Let $X=X_{22}\subset\PP^{13}$ be a Fano threefold of the main series with $\g(X)=12$.
Assume that $\rr(X)> 2$.
Then $X$ contains a surface of degree $d$ with $d\not\equiv 0\mod 11$.
\end{mtheorem}
\begin{proof}
Assume that the degree of any surface $S\subset X$ 
is divisible by $11$ (in particular, $X$ does not contain planes).
Apply the construction \eqref{equation-MMP}.
By Lemma \ref{Lemma-degree-increase}, on each step, 
the variety $\tilde X_k$ does not contain surfaces of degree $\le 10$.
Since the degree of the exceptional divisor of contractions of types \type{e_2}-\type{e_4}
is at most $4$ by \eqref{table-extremal-contractions}, all the birational contractions in 
\eqref{equation-MMP} 
are of type \type{e_1}. We may assume that $\tilde X_N$ does not admit birational contractions.

First consider the case $\uprho(\tilde X_N)=1$.
Then $N\ge 2$ and by 
\eqref{equation-E1-main-0} we have 
$$
64\ge -K_{\tilde X_N}^3\ge 22+20N.
$$
Hence, $N= 2$ and $X_2=\tilde X_2\simeq \PP^3$. 
The morphism $X_1\simeq \tilde X_1\to \tilde X_2\simeq \PP^3$ is 
the blowup of an irreducible curve $B_2\subset \PP^3$.
By \eqref{equation-E1-main}
\begin{equation*}
K_{\tilde X_1}^2\cdot E_1=-K_{\PP^3}\cdot B_2+2-2p_a(B_2)
\end{equation*}
is even.
Hence, $K_{\tilde X_1}^2\cdot E_1\ge 12$ and 
\begin{equation*}
42= (-K_{\PP^3})^3 -(-K_{X_{0}})^3\ge 20+2 K_{\tilde X_1}^2\cdot E_1+2p_a(B_2)-2
\end{equation*}
One can see that the above inequality is in fact an equality and so
$p_a(B_2)=0$, $K_{\tilde X_1}^2\cdot E_1=12$, and
\begin{equation*}
-K_{\PP^3}\cdot B_2=10\not \equiv 0\mod 4, 
\end{equation*}
a contradiction.

Therefore, $\uprho(\tilde X_N)\ge 2$ and $Z$ is not a point.
If 
$\tilde X_N$ has a contraction to a curve, then a general fiber $F$ is a del Pezzo surface. 
Hence, $(-K_{\tilde X})^2\cdot F=K_F^2\le 9$.
This contradicts our assumptions. Thus we may assume that $Z$ is a surface 
with $\uprho(Z)=1$, i.e. $Z\simeq \PP^2$.
Then $\uprho(\tilde X_N)=2$ and so $N\ge 1$.
By \eqref{equation-E1-main-0} 
\begin{equation}\label{equation-KN}
-K_{\tilde X_N}^3\ge 22+ 2\cdot K_{\tilde X}^2\cdot E_0-2\ge 42. 
\end{equation}
Since $\tilde X_N$ has no contractions to a curve and $\uprho(\tilde X_N)\ge 2$, the variety 
$X_N$ is not a quadric. 
By Lemma \ref{Proposition-XN}\ $\iota(\tilde X_N)=2$. 
In this case,
the (anticanonical) degree of any curve on $\tilde X_N$ is even. 
Hence, $(-K_{\tilde X})^2\cdot E_0\equiv 0\mod 22$
by \eqref{equation-E1-main}. This contradicts \eqref{equation-KN}.
\end{proof}

\begin{mtheorem}{\bf Proposition.}\label{Proposition-G-le2}
Let $X=X_{22}\subset\PP^{13}$ be a $G$-Fano threefold of the main series with $\g(X)=12$
and $\rr(X)> 2$.
Let $G_\bullet\subset G$ be an $11$-Sylow subgroup. 
Then $G_\bullet$ acts non-trivially on the lattice $\Cl(X)$. 
In particular, $G_\bullet$ is non-trivial and $\rr(X)\ge 11$.
\end{mtheorem}

\begin{proof}
Let $S_1\subset X$ be a surface of degree $d$.
By Lemma \ref{Lemma-surface-degree=10} we can take $S_1$ so that 
$d\not \equiv 0\mod 11$. 
Let $O= G\cdot S_1$ be the $G$-orbit of $S_1$ and let $n:=\card{O}$.
Denote $D:=\sum_{S\in O} S$. 
We can write $D\sim -aK_X$
for some $a$. Comparing the degrees we get 
\begin{equation*}
nd=22a,\qquad 
n\equiv 0\mod 11.
\end{equation*}
Consider the following natural homomorphism of $G$-modules
\begin{equation*}
\varsigma: \bigoplus_{S\in O} \ZZ\cdot S\longrightarrow \Cl(X).
\end{equation*}
Clearly, $\varsigma(O)$ is a $G$-orbit. Let $\Theta_1,\dots, \Theta_m$ be all elements of
$\varsigma(O)$. Take representatives $S_i\in \varsigma^{-1}(\Theta_i)\cap O$
and put $D':=S_1+\cdots+ S_m$.
Since all the elements in the 
preimage $\varsigma^{-1}(\Theta)\cap O$ are linearly equivalent, 
$m$ divides $n$ and $D\sim \frac nm D'$. 
Hence, 
\begin{equation*}
22D'=\frac {nd}{a}D'\sim \frac {md}aD\sim -mdK_X.
\end{equation*}
Since $\iota(X)=1$, we have $m\equiv 0\mod 11$ and so $G_\bullet$ acts on $\Cl(X)$ non-trivially.
\end{proof}

Now we are in position to prove Theorem \xref{Theorem-main}.

\begin{proof}[Proof of Theorem \xref{Theorem-main}.]
First assume that $\rr(X)= 2$.
Since $\rk \Cl(X)^G=1$, the action of $G$ on $\Cl(X)\simeq \ZZ^2$ is non-trivial. 
Then some element $\tau \in G$ switches two extremal rays $\mathcal R$ and $\mathcal R^+$ of the 
cone of movable divisors $\operatorname{Mov}(X)$.
Hence, two corresponding rational maps $\Phi: X\dashrightarrow Z$ and $\Phi^+: X\dashrightarrow Z^+$
are also switched by $\tau$. In particular, $\dim Z=\dim Z^+$.
This happens only in the case \ref{theorem-v22-rho-2-P3} of Theorem \ref{Theorem-v22-rho-2}.

Now assume that $\rr(X)>2$. 
Then by Proposition \ref{Proposition-G-le2} we have
$\rr(X)\ge 11$. On the other hand, the variety $X$ does not contain planes
(see Theorem \ref{Theorem-planes}) and so we have a contradiction by Proposition
\ref{Proposition-le10}. 
\end{proof}

\input v22-bib.tex

\end{document}

%% file: appendix.tex
\begin{stheorem}{\bf Corollary.}
The Fano threefolds $X_{22}$ of type \xref{theorem-v22-rho-2-DP5} are parametrized by 
the set of $\mathfrak g^1_5$ on $\PP^1$ modulo $\Aut(\PP^1)$.
\end{stheorem}

\begin{stheorem}{\bf Lemma.}\label{Lemma-jumping-lines -v22-IV}
Let $l\subset \PP^2$ be a line.
The following conditions are equivalent:
\begin{enumerate}
\item \label{Lemma-jumping-lines -v22-IV-1}
$l$ is a jumping line of $\EEE$;
\item \label{Lemma-jumping-lines -v22-IV-2}
$l=f(\Lambda)$ for some line $\Lambda\subset Y$;
\item \label{Lemma-jumping-lines -v22-IV-3}
$l$ cuts on $\Omega$ an element $\Xi\in \mathfrak g_5^1$.
\end{enumerate}
\end{stheorem}
\begin{proof}
\ref{Lemma-jumping-lines -v22-IV-1}$\Leftrightarrow$\ref{Lemma-jumping-lines -v22-IV-2}
Let $F:=f^{-1}(l)$ and let $\Lambda$ be the minimal section of $F\simeq\FF_n$.
By the adjunction formula, $K_F=2(M+F)|_F$. Hence $n$ is even. Further,
\[
-K_Y\cdot \Lambda = -K_F\cdot \Lambda+1=3-n.
\]
In the general case we have 
$-K_Y\cdot \Lambda=3$, $n=0$. Thus, the only possibility for the jumping case is the following:
$-K_Y\cdot \Lambda=1$, $n=2$.

The implications \ref{Lemma-jumping-lines -v22-IV-2}$\Leftrightarrow$\ref{Lemma-jumping-lines -v22-IV-3}
are easy and left to the reader.
%
\end{proof}

Note that any line $\Lambda\subset Y$ is contained in a unique member $S_{\Lambda}\in |M+F|$.
\begin{stheorem}{\bf Lemma.}\label{LemmaCorollary-lines-a-v22-IV}
Let $\Lambda\subset Y$ be a line. 
Then $\NNN_{\Lambda/Y}\simeq \OOO_{\PP^1}(-a)\oplus \OOO_{\PP^1}(a-1)$,
where $a=0$ or $-1$. If $\Lambda$ is contained in the smooth locus of $S_{\Lambda}$, then 
$a=0$.
\end{stheorem}

\begin{proof}
The first assertion is standard (see 
e.g. \cite[Lemma 4.2.1]{Iskovskikh-Prokhorov-1999}).
For the second one, note that $\Lambda$ is a $(-1)$-curve on $S_{\Lambda}$
and so $\NNN_{\Lambda/S_{\Lambda}}\simeq \OOO_{\PP^1}(-1)$ is a subbundle of $\NNN_{\Lambda/Y}$. 
\end{proof}

\begin{stheorem}{\bf Corollary.}\label{Corollary-lines-a-v22-IV}
\begin{enumerate}
\item\label{Corollary-lines-a-v22-IV-2}
The scheme $\mathfrak L(Y)$ is purely one-dimensional and generically reduced.
\item\label{Corollary-lines-a-v22-IV-0}
The natural map 
$\mathfrak L(Y)\to (\PP^2)^\vee$, $\Lambda \mapsto f(\Lambda)$ is bijective.

\item\label{Corollary-lines-a-v22-IV-3}
The scheme $\mathfrak L(Y)$ is connected.
\end{enumerate} 
\end{stheorem}
\begin{proof}
For the first (see 
e.g. \cite[Proposition 4.2.2]{Iskovskikh-Prokhorov-1999}).
The rest is easy.
\end{proof}

\begin{scase}
Recall the 
a \emph{L\"uroth quartic} is a plane quartic curve containing the 10 vertices of a complete $5$-side.
All L\"uroth quartics form an irreducible subvariety of codimension one 
in the space of all quartics.
\end{scase}

\begin{stheorem}{\bf Corollary \cite{Barth1977}.}\label{Corollary-jumping-lines-v22}
Let $J(\EEE)\subset {\PP^2}^\vee$ be the curve of jumping lines.
Then $J(\EEE)$ is a L\"uroth quartic curve. 
\end{stheorem}
\begin{proof}
Note that $J(\EEE)$ is the image of the map 
from Corollary \ref{Corollary-lines-a-v22-IV}\ref{Corollary-lines-a-v22-IV-0}.
Take a general member $S\in |M+F|$. Then $S$ is a smooth del Pezzo surface
and the projection $f_S: S\to \PP^2$ contracts $5$ exceptional curves 
$\Upsilon_1,\dots, \Upsilon_5$ (fibers of $f$). 
Let $q_1$, \dots, $q_5\in\Omega\subset \PP^2$ be images of these curves
and let $l_{i,j}$ be the line joining $q_i$ with $q_j$ for $i\neq j$.
The proper transform $\Lambda_{i,j}\subset S$ of $l_{i,j}$
is a $(-1)$-curve meeting $\Upsilon_i$ and $\Upsilon_j$.
Since $\Lambda_{i,j}\cap \Gamma=\emptyset$ and $F\cdot \Lambda_{i,j}=1$,
we have $-K_Y\cdot \Lambda_{i,j}=1$, i.e. $\Lambda_{i,j}$ is a line on $Y$.
The lines $l_{i,j}$ can be regarded as vertices of a complete $5$-side 
on the dual plane ${\PP^2}^\vee$ and the curve of jumping lines 
$J(\EEE)\subset {\PP^2}^\vee$ passes through all these vertices $l_{i,j}$
by Lemma \ref {Lemma-jumping-lines -v22-IV}.

Now take the pencil $l_t$ of lines on $\PP^2$ passing through $q_1$.
Since $\Bs |M+F|=\Gamma$, the surface $S$ is the only member 
of $|M+F|$ containing the fiber $f^{-1}(q_1)$.
Then $l_t$ is a jumping line if and only if there is a line $\Lambda\subset Y$ such that
$f(\Lambda)=l_t$. This $\Lambda$ is contained in some member $S'\in |M+F|$. 
Since $S'\cap f^{-1}(q_1)$ contains two distinct points $\Gamma\cap f^{-1}(q_1)$
and $\Lambda\cap f^{-1}(q_1)$, the surface $S'$ must contain the fiber $f^{-1}(q_1)$.
Thus, $S'=S$ and so $S\supset \Lambda$.
This shows that $\Lambda=\Lambda_{1,j}$ for $j=1,\dots,4$.
Hence, the pencil $l_t$ contains exactly four jumping lines,
i.e. $\deg J(\EEE)=4$.
\end{proof}

%% file: v22-bib.tex
\def\cprime{$'$} \def\mathbb#1{\mathbf#1}